# EXACT LOCAL WHITTLE ESTIMATION OF FRACTIONAL INTEGRATION

By Katsumi Shimotsu[1] and Peter C. B. Phillips[2]

*Queen's University and Yale University, University of Auckland and University of York*

An exact form of the local Whittle likelihood is studied with the intent of developing a general-purpose estimation procedure for the memory parameter ($d$) that does not rely on tapering or differencing prefilters. The resulting exact local Whittle estimator is shown to be consistent and to have the same $N(0, \frac{1}{4})$ limit distribution for all values of $d$ if the optimization covers an interval of width less than $\frac{9}{2}$ and the initial value of the process is known.

**1. Introduction.** Semiparametric estimation of the memory parameter ($d$) in fractionally integrated ($I(d)$) time series is appealing in empirical work because of the general treatment of the short-memory component that it affords. Two common statistical procedures in this class are log-periodogram (LP) regression [1, 10] and local Whittle (LW) estimation [5, 11]. LW estimation is known to be more efficient than LP regression in the stationary ($|d| < \frac{1}{2}$) case, although numerical optimization methods are needed in the calculation. Outside the stationary region, it is known that the asymptotic theory for the LW estimator is discontinuous at $d = \frac{3}{4}$ and again at $d = 1$, is awkward to use because of nonnormal limit theory and, worst of all, the estimator is inconsistent when $d > 1$ [8]. Thus, the LW estimator is not a good general-purpose estimator when the value of $d$ may take on values in the nonstationary zone beyond $\frac{3}{4}$. Similar comments apply in the case of LP estimation [4].

To extend the range of application of these semiparametric methods, data differencing and data tapering have been suggested [3, 15]. These methods

Received March 2002; revised July 2004.
[1]Supported by ESRC Grant R000223629.
[2]Supported by NSF Grant SES-00-92509.
*AMS 2000 subject classification.* 62M10.
*Key words and phrases.* Discrete Fourier transform, fractional integration, long memory, nonstationarity, semiparametric estimation, Whittle likelihood.

have the advantage that they are easy to implement and they make use of existing algorithms once the data filtering has been carried out. Differencing has the disadvantage that prior information is needed on the appropriate order of differencing. Tapering has the disadvantage that the filter distorts the trajectory of the data and inflates the asymptotic variance. As a consequence, there is presently no general-purpose efficient estimation procedure when the value of $d$ may take on values in the nonstationary zone beyond $\frac{3}{4}$.

The present paper studies an exact form of the local Whittle estimator which does not rely on differencing or tapering and which seems to offer a good general-purpose estimation procedure for the memory parameter that applies throughout the stationary and nonstationary regions of $d$. The estimator, which we call the exact LW estimator, is shown to be consistent and to have $N(0, \frac{1}{4})$ limit distribution when the optimization covers an interval of width less than $\frac{9}{2}$. The exact LW estimator therefore has the same limit theory as the LW estimator has for stationary values of $d$. The approach seems to offer a useful alternative for applied researchers who are looking for a general-purpose estimator and want to allow for a substantial range of stationary and nonstationary possibilities for $d$. The method has the further advantage that it provides a basis for constructing asymptotic confidence intervals for $d$ that are valid irrespective of the true value of the memory parameter.

The exact LW estimator given here assumes the initial value of the data to be known. This restriction can be removed by estimating it along with $d$, as shown by Shimotsu [14]. Also, computation of the estimator involves a numerical optimization that is more demanding than conventional LW estimation. Our experience from simulations indicates that the computation time required is about ten times that of the LW estimator and is well within the capabilities of a small notebook computer.

## 2. Exact local Whittle estimation.

We consider the fractional process $X_t$ generated by the model

$$(1 - L)^{d_0} X_t = u_t I\{t \geq 1\}, \qquad t = 0, \pm 1, \ldots, \tag{1}$$

where $I\{\cdot\}$ is the indicator function and $u_t$ is stationary with zero mean and spectral density $f_u(\lambda) \sim G_0$ as $\lambda \to 0$. Expanding the binomial in (1) gives the form

$$\sum_{k=0}^{t} \frac{(-d_0)_k}{k!} X_{t-k} = u_t I\{t \geq 1\}, \tag{2}$$

where

$$(d_0)_k = \frac{\Gamma(d_0 + k)}{\Gamma(d_0)} = (d_0)(d_0 + 1) \cdots (d_0 + k - 1)$$

is Pochhammer's symbol for the forward factorial function and $\Gamma(\cdot)$ is the gamma function. When $d_0$ is a positive integer, the series in (2) terminates, giving the usual formulae for the model (1) in terms of differences and higher-order differences of $X_t$. An alternative form for $X_t$ is obtained by inversion of (1), giving a valid representation for all values of $d_0$,

$$X_t = (1-L)^{-d_0} u_t I\{t \geq 1\} = \sum_{k=0}^{t-1} \frac{(d_0)_k}{k!} u_{t-k}. \tag{3}$$

Define the discrete Fourier transform (d.f.t.) and the periodogram of a time series $a_t$ evaluated at frequency $\lambda$ as

$$w_a(\lambda) = (2\pi n)^{-1/2} \sum_{t=1}^{n} a_t e^{it\lambda},$$
$$I_a(\lambda) = |w_a(\lambda)|^2.$$

2.1. *Exact local Whittle likelihood and estimator.* We start with the likelihood function of the stationary innovation $u_t$. The (negative) Whittle likelihood of $u_t$ based on frequencies up to $\lambda_m$ and up to scale multiplication is

$$\sum_{j=1}^{m} \log f_u(\lambda_j) + \sum_{j=1}^{m} \frac{I_u(\lambda_j)}{f_u(\lambda_j)}, \qquad \lambda_j = \frac{2\pi j}{n}, \tag{4}$$

where $m$ is some integer less than $n$. We want to transform the likelihood function (4) to be data dependent.

If $|d_0| < \frac{1}{2}$, it is known that $I_u(\lambda_j)$ can be approximated by $\lambda_j^{2d_0} I_x(\lambda_j)$ [10, 12]. Therefore, if one views $I_u(\lambda_j)$ as the $j$th observation of $u_t$ in the frequency domain, replacing $I_u(\lambda_j)$ in (4) with $\lambda_j^{2d_0} I_x(\lambda_j)$ and adding the Jacobian $\sum_{j=1}^{m} \log \lambda_j^{-2d}$ to (4) makes it data dependent. Indeed, the resulting objective function coincides with that of the LW estimator.

However, when $d_0$ takes a larger value, in particular when $|d_0| > 1$, $\lambda_j^{2d_0} I_x(\lambda_j)$ no longer provides a good approximation of $I_u(\lambda_j)$. In this paper, we propose to use a "corrected" d.f.t. of $X_t$ that can approximate $I_u(\lambda_j)$ and validly transform (4) in such cases. Lemma 5.1 in Section 5 provides the necessary algebraic relationship for these quantities for any value of $d_0$, namely,

$$\begin{aligned} I_u(\lambda_j) &= I_{\Delta^{d_0} x}(\lambda_j) = |D_n(e^{i\lambda_j}; d_0)|^2 |v_x(\lambda_j; d_0)|^2, \\ v_x(\lambda_j; d) &= w_x(\lambda_j) - D_n(e^{i\lambda_j}; d)^{-1} (2\pi n)^{-1/2} \widetilde{X}_{\lambda_j n}(d), \end{aligned} \tag{5}$$

where

$$D_n(e^{i\lambda}; d) = \sum_{k=0}^{n} \frac{(-d)_k}{k!} e^{ik\lambda}$$

and

$$\widetilde{X}_{\lambda n}(d) = \sum_{p=0}^{n-1} \widetilde{d}_{\lambda p} e^{-ip\lambda} X_{n-p} \qquad \text{with } \widetilde{d}_{\lambda p} = \sum_{k=p+1}^{n} \frac{(-d)_k}{k!} e^{ik\lambda}.$$

The function $v_x(\lambda_j; d_0)$ in (5) adds a correction term that involves $\widetilde{X}_{\lambda_j n}(d_0)$ to the d.f.t. $w_x(\lambda_j)$, which ensures that the relationship (5) holds exactly for all $d_0$. Accordingly, we may interpret $v_x(\lambda_j; d_0)$ as a well-suited proxy for the $j$th frequency domain observation of $X_t$. Consequently, replacing $I_u(\lambda_j)$ in (4) with $|D_n(e^{i\lambda_j}; d)|^2 |v_x(\lambda_j; d)|^2$, adding the Jacobian $\sum_{j=1}^m \log |D_n(e^{i\lambda_j}; d)|^{-2}$ and using (5) again give, in conjunction with the local approximation $f_u(\lambda_j) \sim G$ and $|D_n(e^{i\lambda_j}; d)|^2 \sim \lambda_j^{2d}$ [8],

$$Q_m(G,d) = \frac{1}{m} \sum_{j=1}^{m} \left[ \log(G\lambda_j^{-2d}) + \frac{1}{G} I_{\Delta^d x}(\lambda_j) \right],$$

where $I_{\Delta^d x}(\lambda_j)$ is the periodogram of

$$\Delta^d X_t = (1-L)^d X_t = \sum_{k=0}^{t} \frac{(-d)_k}{k!} X_{t-k}.$$

We propose to estimate $d$ and $G$ by minimizing $Q_m(G,d)$, so that

$$(\widehat{G}, \widehat{d}) = \underset{G \in (0,\infty), d \in [\Delta_1, \Delta_2]}{\arg\min} Q_m(G,d), \tag{6}$$

where $\Delta_1$ and $\Delta_2$ are the lower and upper bounds of the admissible values of $d$ such that $-\infty < \Delta_1 < \Delta_2 < \infty$. Concentrating $Q_m(G,d)$ with respect to $G$, we find that $\widehat{d}$ satisfies

$$\widehat{d} = \underset{d \in [\Delta_1, \Delta_2]}{\arg\min} R(d),$$

where

$$R(d) = \log \widehat{G}(d) - 2d \frac{1}{m} \sum_{j=1}^{m} \log \lambda_j, \qquad \widehat{G}(d) = \frac{1}{m} \sum_{j=1}^{m} I_{\Delta^d x}(\lambda_j).$$

The estimator $\widehat{d}$ is based on the transformation of the Whittle likelihood function of $u_t$ by (5). Since (5) follows from a purely algebraic manipulation and holds exactly for any $d$, we call $\widehat{d}$ the exact local Whittle estimator of $d$. [The word "exact" is used to distinguish the proposed estimator (which relies on an exact algebraic manipulation) from the conventional local Whittle estimator, which is based on the approximation $I_x(\lambda_j) \sim \lambda_j^{-2d} I_u(\lambda_j)$. Of course, the Whittle likelihood is itself an approximation of the exact likelihood, but this should cause no confusion.]

2.2. *Consistency.* We now introduce the assumptions on $m$ and the stationary component $u_t$ in (1).

ASSUMPTION 1.

$$f_u(\lambda) \sim G_0 \in (0, \infty) \qquad \text{as } \lambda \to 0+ .$$

ASSUMPTION 2. In a neighborhood $(0, \delta)$ of the origin, $f_u(\lambda)$ is differentiable and

$$\frac{d}{d\lambda} \log f_u(\lambda) = O(\lambda^{-1}) \qquad \text{as } \lambda \to 0+ .$$

ASSUMPTION 3.

$$u_t = C(L)\varepsilon_t = \sum_{j=0}^{\infty} c_j \varepsilon_{t-j}, \qquad \sum_{j=0}^{\infty} c_j^2 < \infty,$$

where

$$E(\varepsilon_t | F_{t-1}) = 0, \qquad E(\varepsilon_t^2 | F_{t-1}) = 1 \qquad \text{a.s., } t = 0, \pm 1, \ldots,$$

in which $F_t$ is the $\sigma$-field generated by $\varepsilon_s$, $s \le t$, and there exists a random variable $\varepsilon$ such that $E\varepsilon^2 < \infty$ and for all $\eta > 0$ and some $K > 0$, $\Pr(|\varepsilon_t| > \eta) \le K \Pr(|\varepsilon| > \eta)$.

ASSUMPTION 4.

$$\frac{1}{m} + \frac{m(\log m)^{1/2}}{n} + \frac{\log n}{m^\gamma} \to 0 \qquad \text{for any } \gamma > 0.$$

ASSUMPTION 5.

$$\Delta_2 - \Delta_1 \le \tfrac{9}{2}.$$

Assumptions 1–3 are analogous to Assumptions A1–A3 of [11]. However, we impose them in terms of $u_t$ rather than $X_t$. Assumption 4 is slightly stronger than Assumption A4 of [11]. Assumption 5 restricts the length of the interval of permissible values in the optimization (6), although it imposes no restrictions on the value of $d_0$ itself. For instance, if we assume the data are overdifferenced at most once and hence $d_0 \ge -1$, then taking $[\Delta_1, \Delta_2] = [-1, 3.5]$ makes $\widehat{d}$ consistent for any $d_0 \in [\Delta_1, \Delta_2]$. When one wants to allow the interval of permissible values to be wider than $\frac{9}{2}$, the tapered estimators with sufficiently high order of tapering provide useful alternatives.

Under these conditions we may now establish the consistency of $\widehat{d}$.

THEOREM 2.1. *Suppose $X_t$ is generated by (1) with $d_0 \in [\Delta_1, \Delta_2]$ and Assumptions 1–5 hold. Then $\widehat{d} \xrightarrow{p} d_0$ as $n \to \infty$.*

Assumption 5 is necessary for the following reason. Loosely speaking, we prove consistency by showing that (i) when $|d-d_0|$ is small, $R(d)-R(d_0)$ converges uniformly to a nonrandom function that achieves its minimum at $d_0$, and (ii) when $|d-d_0|$ is large, $R(d)-R(d_0)$ is uniformly bounded away from 0. When $|d-d_0|$ is larger than $\frac{1}{2}$, the periodogram $I_{\Delta^d x}(\lambda_j)$ in the objective function does not behave like $\lambda_j^{2(d-d_0)} I_u(\lambda_j)$. Consequently, $R(d)-R(d_0)$ does not converge to a nonrandom function, and we need an alternative way to bound it away from zero. For instance, when $\frac{1}{2} \le d-d_0 \le \frac{3}{2}$, the normalized d.f.t. is expressed as [cf. (30) in the proof of consistency]

$$\lambda_j^{-(d-d_0)} w_{\Delta^d x}(\lambda_j) \simeq e^{-(\pi/2)(d-d_0)i} w_u(\lambda_j) + \lambda_j^{-(d-d_0)} (2\pi n)^{-1/2} e^{i\lambda_j} Z_n,$$

where

$$Z_n = \sum_{t=1}^{n} (1-L)^d X_t.$$

The leakage from the last term prevents the uniform convergence of $R(d)-R(d_0)$ and complicates the proof. When $|d-d_0|$ is larger, $\lambda_j^{-(d-d_0)} w_{\Delta^d x}(\lambda_j)$ has further additional terms [e.g., the equation below (51)], and we were able to show the necessary results only for $|d-d_0| \le \frac{9}{2}$, which is why we need Assumption 5. Lemma 5.10 in Section 5 is the main tool in handling the effects of such additional terms. We could relax Assumption 5 if we could extend Lemma 5.10 to hold with more general summands $(1-e^{i\lambda_j})^k Q_k + \cdots + Q_0$, but we were not able to do so.

Remark 1. An alternative way of accommodating a wider range of $d$ without sacrificing efficiency is to use a two-step procedure. A two-step estimator based on the objective function $R(d)$ that uses a (higher-order) tapered estimator in the first step would have the same asymptotic variance as the exact LW estimator. (Strictly speaking, the asymptotic properties of tapered estimators have been established only under the alternative type of fractionally integrated process generated as in (8), although some results on the difference between their d.f.t.'s are available [12].)

Remark 2. The model (1) assumes, in effect, that the initial value of $X_t$ is known. In practice, it is more natural to allow for an unknown initial value, $\mu_0$, and model $X_t$ as

$$X_t = \mu_0 + (1-L)^{-d_0} u_t I\{t \ge 1\} = \mu_0 + \sum_{k=0}^{t-1} \frac{(d_0)_k}{k!} u_{t-k}. \tag{7}$$

Estimation of $\mu_0$ affects the limiting behavior of the estimator. According to Shimotsu [14], (i) if $\mu_0$ is replaced with the sample average $\bar{X} = n^{-1}\sum_{t=1}^n X_t$, then the estimator is consistent for $d_0 \in (-\frac{1}{2}, 1)$ and asymptotically normal for $d_0 \in (-\frac{1}{2}, \frac{3}{4})$, but simulations suggest that the estimator is inconsistent for $d_0 > 1$; and (ii) if $\mu_0$ is replaced by $X_1$, then the estimator is consistent for $d_0 \geq \frac{1}{2}$ and asymptotically normal for $d_0 \in [\frac{1}{2}, 2)$, but simulations suggest that the estimator is inconsistent for $d_0 \leq 0$. To accommodate unknown $\mu_0$, it is possible to extend Theorem 2.1 for $X_t$ generated by (7) by estimating $\mu_0$ along with $d_0$. For instance, Shimotsu [14] proposes estimating $\mu_0$ by

$$\widehat{\mu}(d) = w(d)\bar{X} + (1 - w(d))X_1,$$

where $w(d)$ is a smooth (twice continuously differentiable) weight function such that $w(d) = 1$ for $d \leq \frac{1}{2}$, $w(d) \in [0, 1]$ for $\frac{1}{2} \leq d \leq \frac{3}{4}$ and $w(d) = 0$ for $d \geq \frac{3}{4}$, and replacing $X_t$ with $X_t - \widehat{\mu}(d)$ in the periodograms in the objective function. Shimotsu [14] shows the resulting estimator of $d$ is consistent and asymptotically normal for $d_0 \in (-\frac{1}{2}, 2)$, excluding arbitrary small intervals around 0 and 1. Another possibility would be to replace $X_t$ with $X_t - \mu$ in the periodogram ordinates and minimize the objective function with respect to $(d, G, \mu)$.

Remark 3. Fractionally integrated processes as defined in (1) are more restrictive in some ways than the stationary frequency domain characterization used in [11] and elsewhere. It might be possible to extend the results in this paper to the class of nonstationary processes analyzed by [13] and seek to achieve a similar degree of generality to Robinson [11], but we do not attempt to do so here.

Remark 4. Another popular definition of a fractionally integrated process provides for different generating mechanisms according to the specific range of values taken by $d_0$, as in

$$X_t = \begin{cases} (1-L)^{-d_0} u_t, & d_0 \in (-\infty, \frac{1}{2}), \\ \mu_0 + \sum_{k=1}^{t} Z_k, \qquad Z_t = (1-L)^{1-d_0} u_t, & d_0 \in [\frac{1}{2}, \frac{3}{2}), \end{cases} \tag{8}$$

with corresponding extensions for larger values of $d_0$, so that $X_t$ or its (higher-order) difference is stationary. While we do not explore the effects of these alternative generating mechanisms here, simulation results suggest that the version of the exact LW estimator in [14] is consistent for this type of fractionally integrated process.

2.3. *Asymptotic normality.* We introduce some further assumptions that are used to derive the limit distribution theory.

ASSUMPTION 1′. Assumption 1 holds, and also for some $\beta \in (0,2]$

$$f_u(\lambda) = G_0(1 + O(\lambda^\beta)) \qquad \text{as } \lambda \to 0+ .$$

ASSUMPTION 2′. In a neighborhood $(0,\delta)$ of the origin, $C(e^{i\lambda})$ is differentiable and

$$\frac{d}{d\lambda} C(e^{i\lambda}) = O(\lambda^{-1}) \qquad \text{as } \lambda \to 0+ .$$

ASSUMPTION 3′. Assumption 3 holds and also

$$E(\varepsilon_t^3 | F_{t-1}) = \mu_3, \qquad E(\varepsilon_t^4 | F_{t-1}) = \mu_4 \qquad \text{a.s., } t = 0, \pm 1, \ldots,$$

for finite constants $\mu_3$ and $\mu_4$.

ASSUMPTION 4′. As $n \to \infty$,

$$\frac{1}{m} + \frac{m^{1+2\beta}(\log m)^2}{n^{2\beta}} + \frac{\log n}{m^\gamma} \to 0 \qquad \text{for any } \gamma > 0.$$

ASSUMPTION 5′. Assumption 5 holds.

Assumptions 1′–3′ are analogous to Assumptions A1′–A3′ of [11], except that our assumptions are in terms of $u_t$ rather than $X_t$. Assumption 4′ is slightly stronger than Assumption 4′ of [11].

The following theorem establishes the asymptotic normality of the exact local Whittle estimator for $d_0 \in (\Delta_1, \Delta_2)$. (The approximate mean squared error and the corresponding optimal bandwidth can be obtained heuristically in the same manner as in [2].)

THEOREM 2.2. *Suppose $X_t$ is generated by (1) with $d_0 \in (\Delta_1, \Delta_2)$ and Assumptions 1′–5′ hold. Then*

$$m^{1/2}(\widehat{d} - d_0) \xrightarrow{d} N(0, \tfrac{1}{4}) \qquad \text{as } n \to \infty.$$

**3. Simulations.** This section reports some simulations that were conducted to examine the finite sample performance of the exact LW estimator (hereafter, exact estimator), the LW estimator (hereafter, untapered estimator) and the LW estimator with two types of tapering studied by Hurvich and Chen [3] and Velasco [15] with Bartlett's window [hereafter, tapered (HC) and tapered (V) estimator, resp.]. The tapered (HC) estimator and tapered (V) estimator are consistent and asymptotically normal for $d \in (-0.5, 1.5)$ with limiting variances $1.5/(4m)$ and $2.1/(4m)$, respectively (see Remark 1). We generate $I(d)$ processes according to (3) with $u_t \sim \text{i.i.d.} N(0,1)$. $\Delta_1$ and $\Delta_2$ are set to $-6$ and 6. Although this setting violates Assumption 5, it does not appear to adversely affect the performance of the exact estimator. The bias, standard deviation and mean squared error (MSE) were computed using 10,000 replications. The sample size and band parameter $m$ were chosen to be $n = 500$ and $m = n^{0.65} = 56$. Values of $d$ were selected in the interval $[-3.5, 3.5]$.

TABLE 1
*Simulation results:* $n = 500$, $m = n^{0.65} = 56$

| | Exact estimator | | | Untapered estimator | | |
|---|---|---|---|---|---|---|
| **d** | **bias** | **s.d.** | **MSE** | **bias** | **s.d.** | **MSE** |
| −3.5 | −0.0024 | 0.0787 | 0.0062 | 3.1617 | 0.2831 | 10.076 |
| −2.3 | −0.0020 | 0.0774 | 0.0060 | 1.6345 | 0.3041 | 2.7640 |
| −1.7 | −0.0020 | 0.0776 | 0.0060 | 0.8709 | 0.2788 | 0.8363 |
| −1.3 | −0.0014 | 0.0770 | 0.0059 | 0.4109 | 0.2170 | 0.2160 |
| −0.7 | −0.0024 | 0.0787 | 0.0062 | 0.0353 | 0.0885 | 0.0091 |
| −0.3 | −0.0033 | 0.0777 | 0.0060 | −0.0027 | 0.0781 | 0.0061 |
| 0.0 | −0.0029 | 0.0784 | 0.0061 | −0.0075 | 0.0781 | 0.0062 |
| 0.3 | −0.0020 | 0.0782 | 0.0061 | −0.0066 | 0.0785 | 0.0062 |
| 0.7 | −0.0017 | 0.0777 | 0.0060 | 0.0099 | 0.0812 | 0.0067 |
| 1.3 | −0.0014 | 0.0781 | 0.0061 | −0.2108 | 0.0982 | 0.0541 |
| 1.7 | −0.0025 | 0.0780 | 0.0061 | −0.6288 | 0.1331 | 0.4130 |
| 2.3 | −0.0026 | 0.0772 | 0.0060 | −1.2647 | 0.1046 | 1.6104 |
| 3.5 | −0.0016 | 0.0770 | 0.0059 | −2.4919 | 0.0724 | 6.2150 |

Tables 1 and 2 show the simulation results. The exact estimator has little bias for all values of $d$. The untapered estimator has a large bias for $d > 1$, corroborating the theoretical result that it converges to unity in probability [8]. When $-0.5 < d < 1$, the exact and untapered estimators have similar variance and MSE. The variances of the tapered estimators are always larger than those of the exact estimator. Again, this outcome corroborates the theoretical result that the tapered estimators have larger asymptotic variance. The tapered (HC) estimator has small bias and performs better than the tapered (V) estimator for $-0.5 < d < 2$. However, the tapered (HC) estimator has around 50% larger MSE than the exact estimator even for those values of $d$ due to its large variance.

Figures 1 and 2 plot kernel estimates of the densities of the four estimators for the values $d = -0.7, 0.3, 1.3$ and 2.3. The sample size and $m$ were chosen as $n = 500$ and $m = n^{0.65}$, and 10,000 replications were used. When $d = -0.7$, the exact and tapered (V) estimators have symmetric distributions centered on $-0.7$, with the tapered estimator having a flatter distribution. The untapered and tapered (HC) estimators appear to be biased. When $d = 0.3$, the untapered and exact estimators have almost identical distributions, whereas the two tapered estimators have more dispersed distributions. When $d = 1.3$, the untapered estimator is centered on unity. In this case, the exact estimator seems to work well, having a symmetric distribution centered on 1.3. The tapered estimators have flatter distributions than the exact estimator but otherwise appear reasonable and they are certainly better than the inconsistent untapered estimator. When $d = 2.3$, the untapered and tapered (V) estimators appear centered on 1.0 and 2.0, respectively. In this case, the tapered (HC) estimator is upward biased. Again, the exact estimator has a symmetric distribution centered on the true value 2.3.

In summary, there seems to be little doubt from these results that the exact LW estimator is the best general-purpose estimator over a wide range of $d$ values.

## 4. Proofs.

In this and the following section, $|x|_+$ denotes $\max\{x, 1\}$ and $x^*$ denotes the complex conjugate of $x$. $C$, $c$ and $\varepsilon$ denote generic constants such that $C, c \in (1, \infty)$ and $\varepsilon \in (0, 1)$ unless specified otherwise, and they may take different values in different places.

4.1. *Proof of consistency.* Define $G(d) = G_0 \frac{1}{m} \sum_1^m \lambda_j^{2(d-d_0)}$ and $S(d) = R(d) - R(d_0)$. Rewrite $S(d)$ as

$$S(d) = R(d) - R(d_0)$$

TABLE 2
*Simulation results: $n = 500$, $m = n^{0.65} = 56$*

| | **Tapered (HC) estimator** | | | **Tapered (V) estimator** | | |
|---|---|---|---|---|---|---|
| ***d*** | **bias** | **s.d.** | **MSE** | **bias** | **s.d.** | **MSE** |
| −3.5 | 2.5889 | 0.3037 | 6.7946 | 1.6126 | 0.3380 | 2.7148 |
| −2.3 | 1.1100 | 0.2893 | 1.3157 | 0.2155 | 0.1726 | 0.0762 |
| −1.7 | 0.4474 | 0.2154 | 0.2466 | 0.0259 | 0.1235 | 0.0159 |
| −1.3 | 0.1551 | 0.1231 | 0.0392 | 0.0081 | 0.1211 | 0.0147 |
| −0.7 | 0.0278 | 0.0957 | 0.0099 | −0.0068 | 0.1219 | 0.0149 |
| −0.3 | 0.0100 | 0.0971 | 0.0095 | −0.0133 | 0.1224 | 0.0151 |
| 0.0 | 0.0034 | 0.0985 | 0.0097 | −0.0138 | 0.1224 | 0.0152 |
| 0.3 | −0.0033 | 0.1004 | 0.0101 | −0.0132 | 0.1235 | 0.0154 |
| 0.7 | −0.0066 | 0.0994 | 0.0099 | −0.0068 | 0.1227 | 0.0151 |
| 1.3 | −0.0079 | 0.0987 | 0.0098 | 0.0140 | 0.1232 | 0.0154 |
| 1.7 | 0.0008 | 0.0972 | 0.0095 | 0.0456 | 0.1288 | 0.0187 |
| 2.3 | 0.0528 | 0.0981 | 0.0124 | −0.1781 | 0.1419 | 0.0519 |
| 3.5 | −0.4079 | 0.1142 | 0.1795 | −1.4541 | 0.1338 | 2.1322 |

$$= \log \frac{\widehat{G}(d)}{G(d)} - \log \frac{\widehat{G}(d_0)}{G_0} + \log\left(\frac{1}{m}\sum_{j=1}^{m} j^{2d-2d_0} \Big/ \frac{m^{2(d-d_0)}}{2(d-d_0)+1}\right)$$

$$- 2(d-d_0)\left[\frac{1}{m}\sum_{j=1}^{m}\log j - (\log m - 1)\right]$$

$$+ 2(d-d_0) - \log(2(d-d_0)+1).$$

Define $U(d) = 2(d-d_0) - \log(2(d-d_0)+1)$ and

$$T(d) = \log \frac{\widehat{G}(d_0)}{G_0} - \log \frac{\widehat{G}(d)}{G(d)} - \log\left(\frac{1}{m}\sum_{j=1}^{m} j^{2d-2d_0} \Big/ \frac{m^{2d-2d_0}}{2(d-d_0)+1}\right)$$

$$+ 2(d-d_0)\left[\frac{1}{m}\sum_{j=1}^{m}\log j - (\log m - 1)\right],$$

so that $S(d) = U(d) - T(d)$. For arbitrarily small $\Delta > 0$, define $\Theta_1 = \{d_0 - \frac{1}{2} + \Delta \le d \le d_0 + \frac{1}{2}\}$ and $\Theta_2 = \{d \in [\Delta_1, d_0 - \frac{1}{2} + \Delta] \cup [d_0 + \frac{1}{2}, \Delta_2]\}$, $\Theta_2$ being possibly empty. Without loss of generality we assume $\Delta < \frac{1}{8}$ hereafter. For $\frac{1}{2} > \rho > 0$, define $N_\rho = \{d : |d - d_0| < \rho\}$. Then it follows (cf. [11], page 1634) that

$$\Pr(|\widehat{d} - d_0| \ge \rho) \le \Pr\Big(\inf_{d \in \Theta_1 \setminus N_\rho} S(d) \le 0\Big) + \Pr\Big(\inf_{\Theta_2} S(d) \le 0\Big). \tag{9}$$

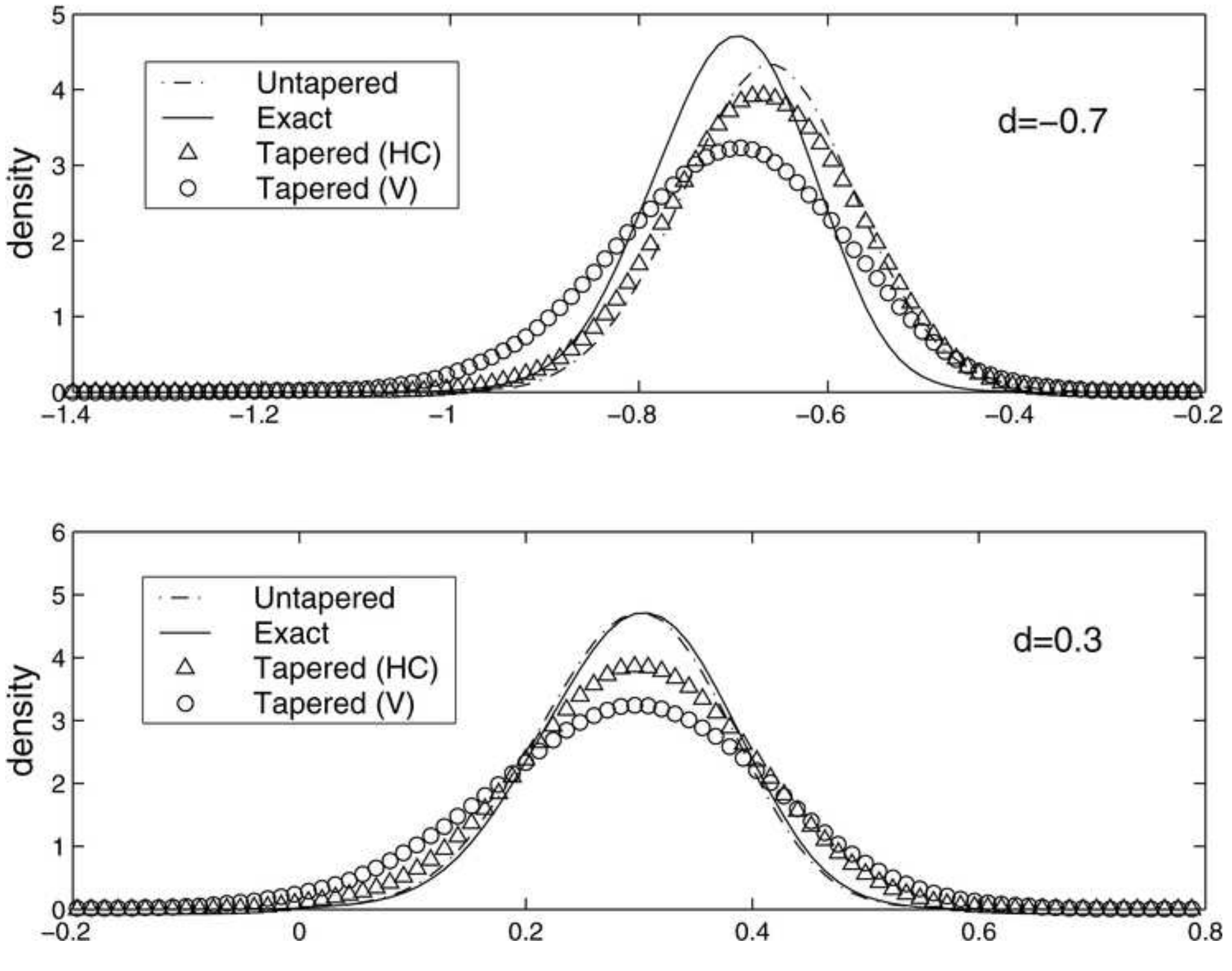

FIG. 1. *Densities of the four estimators:* $n = 500$, $m = n^{0.65}$.

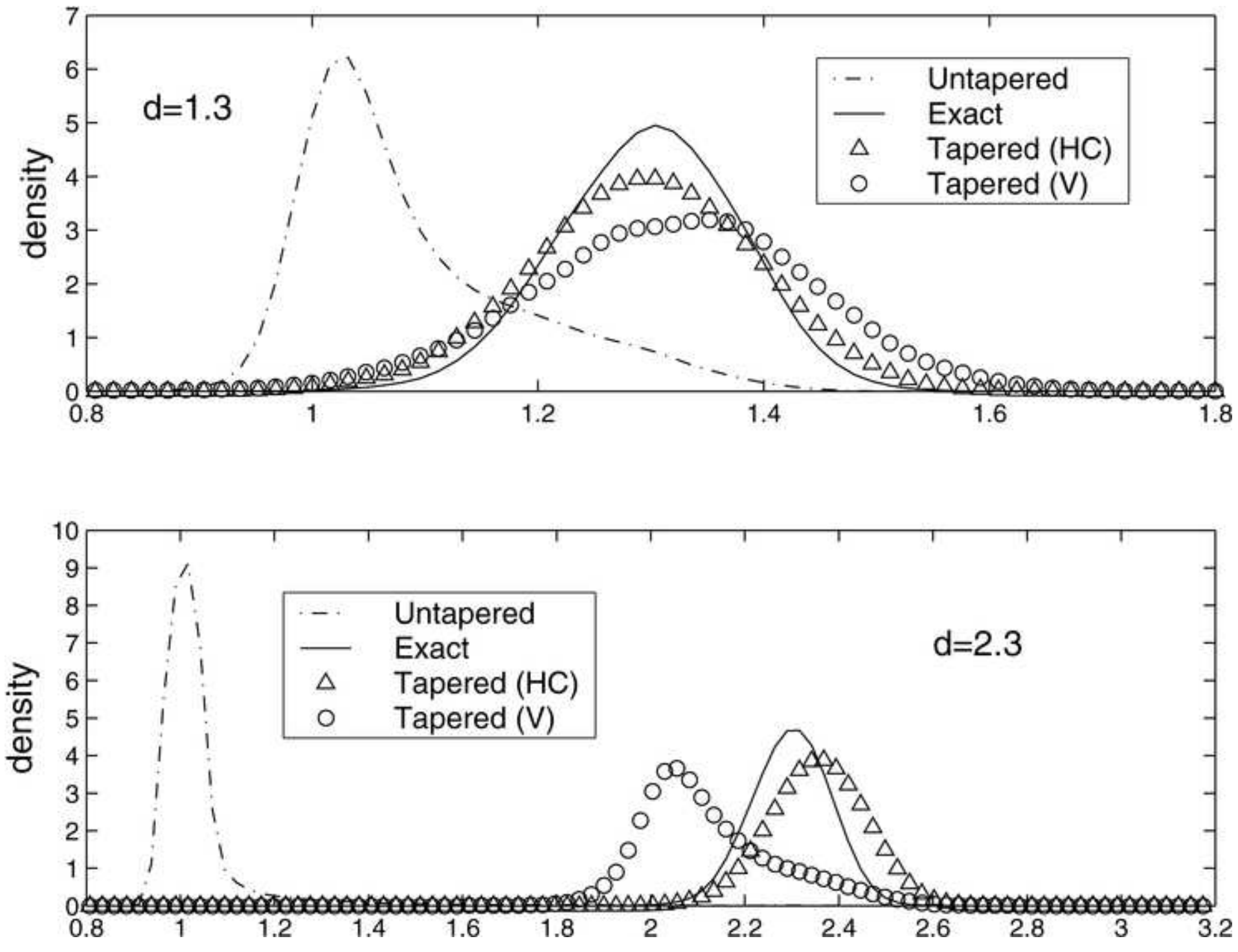

FIG. 2. *Densities of the four estimators:* $n = 500$, $m = n^{0.65}$.

Robinson ([11], (3.4), page 1635) shows

$$\inf_{d \in \Theta_1 \setminus N_\rho} U(d) \ge \rho^2/2. \tag{10}$$

Therefore, $\Pr(|\widehat{d} - d_0| \ge \rho) \to 0$ if

$$\sup_{\Theta_1} |T(d)| \xrightarrow{p} 0, \qquad \Pr\Big(\inf_{\Theta_2} S(d) \le 0\Big) \to 0$$

as $n \to \infty$. From [11], the fourth term of $T(d)$ is $O(\log m/m)$ uniformly in $d \in \Theta_1$ and

$$\sup_{\Theta_1} \left| \frac{2(d - d_0) + 1}{m} \sum_{j=1}^{m} \left(\frac{j}{m}\right)^{2d - 2d_0} - 1 \right| = O\left(\frac{1}{m^{2\Delta}}\right). \tag{11}$$

Note that

$$
\begin{aligned}
&\frac{\widehat{G}(d)-G(d)}{G(d)} \\
&\quad = \frac{m^{-1}\sum_1^m \lambda_j^{2(d-d_0)}\lambda_j^{2(d_0-d)} I_{\Delta^d x}(\lambda_j) - G_0 m^{-1}\sum_1^m \lambda_j^{2(d-d_0)}}{G_0 m^{-1}\sum_1^m \lambda_j^{2(d-d_0)}} \\
&\quad = \frac{m^{-1}\sum_1^m (j/m)^{2(d-d_0)}\lambda_j^{2(d_0-d)} I_{\Delta^d x}(\lambda_j) - G_0 m^{-1}\sum_1^m (j/m)^{2(d-d_0)}}{G_0 m^{-1}\sum_1^m (j/m)^{2(d-d_0)}} \qquad (12)\\
&\quad = \frac{[2(d-d_0)+1]m^{-1}\sum_1^m (j/m)^{2(d-d_0)}[\lambda_j^{2(d_0-d)} I_{\Delta^d x}(\lambda_j) - G_0]}{[2(d-d_0)+1]G_0 m^{-1}\sum_1^m (j/m)^{2(d-d_0)}} \\
&\quad = \frac{A(d)}{B(d)}.
\end{aligned}
$$

Therefore, by the fact that $\Pr(|\log Y| \geq \varepsilon) \leq 2\Pr(|Y-1| \geq \varepsilon/2)$ for any nonnegative random variable $Y$ and $\varepsilon \leq 1$, $\sup_{\Theta_1}|T(d)| \xrightarrow{p} 0$ if

$$
\sup_{\Theta_1} |A(d)/B(d)| \xrightarrow{p} 0. \qquad (13)
$$

Define $\theta = d - d_0$, and define

$$
Y_t(\theta) = (1-L)^d X_t = (1-L)^{d-d_0}(1-L)^{d_0} X_t = (1-L)^{\theta} u_t I\{t \geq 1\}.
$$

Hereafter, we use the notation $a_t \sim I(\alpha)$ when $a_t$ is generated by (1) with parameter $\alpha$. So $Y_t \sim I(-\theta)$. Note that

$$
d \in \Theta_1 \quad \Longleftrightarrow \quad -\tfrac{1}{2} + \Delta \leq \theta \leq \tfrac{1}{2}.
$$

Applying Lemma 5.1(a) to $(Y_t(\theta), u_t)$ and reversing the role of $X_t$ and $u_t$, we obtain

$$
w_y(\lambda_j) = w_u(\lambda_j) D_n(e^{i\lambda_j};\theta) - (2\pi n)^{-1/2}\widetilde{U}_{\lambda_j n}(\theta), \qquad (14)
$$

and $A(d)$ can be written as, with $g = 2(d-d_0)+1$,

$$
A(d) = \frac{g}{m}\sum_{j=1}^m \left(\frac{j}{m}\right)^{2\theta}[\lambda_j^{-2\theta} I_y(\lambda_j) - G_0]. \qquad (15)
$$

Hereafter in this section let $I_{yj}$ denote $I_y(\lambda_j)$, let $w_{uj}$ denote $w_u(\lambda_j)$, and employ the same notation for the other d.f.t.'s and periodograms. From an argument similar to that of [11], page 1636, $\sup_{\Theta_1}|A(d)|$ is bounded by

$$
2\sum_{r=1}^{m-1}\left(\frac{r}{m}\right)^{2\Delta}\frac{1}{r^2}\sup_{\Theta_1}\left|\sum_{j=1}^r [\lambda_j^{-2\theta} I_{yj} - G_0]\right| + \frac{2}{m}\sup_{\Theta_1}\left|\sum_{j=1}^m [\lambda_j^{-2\theta} I_{yj} - G_0]\right|. \qquad (16)
$$

Now

$$\begin{aligned} &\lambda_j^{-2\theta} I_{yj} - G_0 \\ &\quad = \lambda_j^{-2\theta} I_{yj} - \lambda_j^{-2\theta} |D_n(e^{i\lambda_j};\theta)|^2 I_{uj} \\ &\qquad + [\lambda_j^{-2\theta} |D_n(e^{i\lambda_j};\theta)|^2 - G_0/f_u(\lambda_j)] I_{uj} \\ &\qquad + [I_{uj} - |C(e^{i\lambda_j})|^2 I_{\varepsilon j}] G_0/f_u(\lambda_j) + G_0(2\pi I_{\varepsilon j} - 1). \end{aligned} \tag{17}$$

For any $\eta > 0$, Lemma 5.2 and Assumption 1 imply that $n$ can be chosen so that

$$|\lambda_j^{-2\theta} |D_n(e^{i\lambda_j};\theta)|^2 - G_0/f_u(\lambda_j)| \le \eta + O(\lambda_j^2) + O(j^{-1/2}), \qquad j = 1, \dots, m. \tag{18}$$

The results in [11], page 1637, imply that, uniformly in $j = 1, \dots, m$,

$$\begin{aligned} &E|w_{uj} - C(e^{i\lambda_j}) w_{\varepsilon j}|^2 = O(j^{-1} \log(j+1)), \\ &E|I_{uj} - |C(e^{i\lambda_j})|^2 I_{\varepsilon j}| = O(j^{-1/2} (\log(j+1))^{1/2}). \end{aligned} \tag{19}$$

It follows from (18) and (19) that

$$\begin{aligned} \sum_{r=1}^m \left(\frac{r}{m}\right)^{2\Delta} \frac{1}{r^2} \sup_{\Theta_1} \sum_{j=1}^r &|[\lambda_j^{-2\theta} |D_n(e^{i\lambda_j};\theta)|^2 - G_0/f_u(\lambda_j)] I_{uj} \\ &+ [I_{uj} - |C(e^{i\lambda_j})|^2 I_{\varepsilon j}] G_0/f_u(\lambda_j)| \\ &= O_p(\eta + m^2 n^{-2} + m^{-2\Delta} \log m). \end{aligned}$$

Robinson ([11], pages 1637–1638) shows $\sum_1^m (r/m)^{2\Delta} r^{-2} |\sum_1^r (2\pi I_{\varepsilon j} - 1)| \xrightarrow{p} 0$ and $m^{-1} \sum_1^m (2\pi I_{\varepsilon j} - 1) \xrightarrow{p} 0$. From (14), the fact that $||A|^2 - |B|^2| \le |A + B||A - B|$ and the Cauchy–Schwarz inequality we have

$$\begin{aligned} &E \sup_{\Theta_1} |\lambda_j^{-2\theta} I_{yj} - \lambda_j^{-2\theta} |D_n(e^{i\lambda_j};\theta)|^2 I_{uj}| \\ &\quad \le \left( E \sup_{\Theta_1} \left| 2\lambda_j^{-\theta} D_n(e^{i\lambda_j};\theta) w_{uj} - \lambda_j^{-\theta} \frac{\widetilde{U}_{\lambda_j n}(\theta)}{\sqrt{2\pi n}} \right|^2 \right)^{1/2} \\ &\qquad \times \left( E \sup_{\Theta_1} \left| \lambda_j^{-\theta} \frac{\widetilde{U}_{\lambda_j n}(\theta)}{\sqrt{2\pi n}} \right|^2 \right)^{1/2}. \end{aligned} \tag{20}$$

From (19) and Lemmas 5.2 and 5.3, it follows that, uniformly in $j = 1, \dots, m$,

$$E \sup_{\Theta_1} |\lambda_j^{-\theta} D_n(e^{i\lambda_j};\theta) w_{uj}|^2 = O(1),$$

$$E \sup_{\Theta_1} |\lambda_j^{-\theta} (2\pi n)^{-1/2} \widetilde{U}_{\lambda_j n}(\theta)|^2 = O(j^{-1} (\log n)^2).$$

Therefore, we obtain

$$(20) = O(1 + j^{-1/2}\log n)O(j^{-1/2}\log n) = O(j^{-1/2}(\log n)^2). \tag{21}$$

It follows that

$$\begin{aligned}\sum_{r=1}^{m-1}\left(\frac{r}{m}\right)^{2\Delta}\frac{1}{r^2}E\sup_{\Theta_1}\left|\sum_{j=1}^{r}[\lambda_j^{-2\theta}I_{yj} - \lambda_j^{-2\theta}|D_n(e^{i\lambda_j};\theta)|^2 I_{uj}]\right| \\ = O(m^{-\Delta}(\log n)^2),\end{aligned}$$

and hence the first term in (16) is $o_p(1)$. Using the same technique, we can show that the second term in (16) is $o_p(1)$, and $\sup_{\Theta_1}|A(d)| \xrightarrow{p} 0$ follows. Equation (11) gives $\sup_{\Theta_1}|B(d) - G_0| = O(m^{-2\Delta})$, and (13) follows.

Now we take care of $\Theta_2 = \{d \in [\Delta_1, d_0 - \frac{1}{2} + \Delta] \cup [d_0 + \frac{1}{2}, \Delta_2]\} = \{\theta \in [\Delta_1 - d_0, -\frac{1}{2} + \Delta] \cup [\frac{1}{2}, \Delta_2 - d_0]\}$ to show $\Pr(\inf_{\Theta_2} S(d) \le 0) \to 0$. Note that

$$\begin{aligned}S(d) &= \log\widehat{G}(d) - \log\widehat{G}(d_0) - 2(d - d_0)\frac{1}{m}\sum_{j=1}^{m}\log\lambda_j \\ &= \log\frac{1}{m}\sum_{j=1}^{m}I_{\Delta^d xj} - \log\frac{1}{m}\sum_{j=1}^{m}I_{\Delta^{d_0}xj} \\ &\quad - 2(d - d_0)\log\frac{2\pi}{n} - 2(d - d_0)\frac{1}{m}\sum_{j=1}^{m}\log j \\ &= \log\frac{1}{m}\sum_{j=1}^{m}\lambda_j^{2(d-d_0)}\lambda_j^{2(d_0-d)}I_{\Delta^d xj} - \log\frac{1}{m}\sum_{j=1}^{m}I_{\Delta^{d_0}xj} \\ &\quad - 2(d - d_0)\log\frac{2\pi}{n} - 2(d - d_0)\log p \\ &= \log\frac{1}{m}\sum_{j=1}^{m}\left(\frac{j}{p}\right)^{2\theta}\lambda_j^{-2\theta}I_{\Delta^d xj} - \log\frac{1}{m}\sum_{j=1}^{m}I_{\Delta^{d_0}xj} \\ &= \log\widehat{D}(d) - \log\widehat{D}(d_0),\end{aligned}$$

where $p = \exp(m^{-1}\sum_1^m \log j) \sim m/e$ as $m \to \infty$. Applying (17) with $\theta = 0$ and proceeding similarly to the argument below (17), we obtain

$$\log\widehat{D}(d_0) - \log G_0 = \log\left(1 + G_0^{-1}\left(\frac{1}{m}\sum_{j=1}^{m}I_{uj} - G_0\right)\right) = o_p(1).$$

Therefore, $\Pr(\inf_{\Theta_2} S(d) \le 0)$ tends to 0 if there exists $\delta > 0$ such that

$$\Pr\left(\inf_{\Theta_2}\log\widehat{D}(d) - \log G_0 \le \log(1+\delta)\right) = \Pr\left(\inf_{\Theta_2}\widehat{D}(d) - G_0 \le \delta G_0\right) \to 0$$

as $n \to 0$. Now, for any fixed $\kappa \in (0,1)$ we have

$$\widehat{D}(d) = \frac{1}{m}\sum_{j=1}^{m}\left(\frac{j}{p}\right)^{2\theta}\lambda_j^{-2\theta}I_{yj} \geq \frac{1}{m}\sum_{j=[\kappa m]}^{m}\left(\frac{j}{p}\right)^{2\theta}\lambda_j^{-2\theta}I_{yj}.$$

Let $\sum'$ denote the sum over $j = [\kappa m], \ldots, m$. It follows that, for $d \in \Theta_2$,

$$\begin{aligned}&\widehat{D}(d) - G_0\\ &\quad \geq m^{-1}{\sum}'(j/p)^{2\theta}(\lambda_j^{-2\theta}I_{yj} - G_0) + G_0\left(m^{-1}{\sum}'(j/p)^{2\theta} - 1\right).\end{aligned} \tag{22}$$

From Lemma 5.5, by choosing $\delta$ first and then $\kappa$ sufficiently small, for large $m$ we have

$$\inf_{\Theta_2} G_0\left(m^{-1}{\sum}'(j/p)^{2\theta} - 1\right) > 4\delta G_0.$$

Therefore, $\Pr(\inf_{\Theta_2} S(d) \leq 0) \to 0$ if there exists $\delta > 0$ such that

$$\Pr\left(\inf_{\Theta_2}\left(m^{-1}{\sum}'(j/p)^{2\theta}(\lambda_j^{-2\theta}I_{yj} - G_0)\right) \leq -3\delta G_0\right) \to 0 \tag{23}$$

as $n \to \infty$. We proceed to show (23) for subsets of $\Theta_2$.

First we consider $\Theta_2^1 = \{\theta \in [-\frac{1}{2}, -\frac{1}{2} + \Delta]\}$. Rewrite

$$m^{-1}{\sum}'(j/p)^{2\theta}(\lambda_j^{-2\theta}I_{yj} - G_0) = \Lambda_{1n}(\theta) + \Lambda_{2n}(\theta),$$

where

$$\Lambda_{1n}(\theta) = m^{-1}{\sum}'(j/p)^{2\theta}[\lambda_j^{-2\theta}I_{yj} - \lambda_j^{-2\theta}|D_n(e^{i\lambda_j};\theta)|^2 I_{uj}], \tag{24}$$

$$\Lambda_{2n}(\theta) = m^{-1}{\sum}'(j/p)^{2\theta}[\lambda_j^{-2\theta}|D_n(e^{i\lambda_j};\theta)|^2 I_{uj} - G_0]. \tag{25}$$

For $\Lambda_{1n}(\theta)$, since (20) and (21) are still valid for $\theta \in \Theta_2^1$, we have

$$E\sup_{\Theta_2^1}|\lambda_j^{-2\theta}I_{yj} - \lambda_j^{-2\theta}|D_n(e^{i\lambda_j};\theta)|^2 I_{uj}| = O(j^{-1/2}(\log n)^2),$$

and it follows from Lemma 5.4 that $E\sup_{\Theta_2^1}|\Lambda_{1n}(\theta)| = o(1)$. For $\Lambda_{2n}(\theta)$, rewrite $\Lambda_{2n}(\theta)$ as

$$m^{-1}{\sum}'(j/p)^{2\theta}[\lambda_j^{-2\theta}|D_n(e^{i\lambda_j};\theta)|^2 - G_0/f_u(\lambda_j)]I_{uj} \tag{26}$$

$$+ m^{-1}{\sum}'(j/p)^{2\theta}[I_{uj} - |C(e^{i\lambda_j})|^2 I_{\varepsilon j}]G_0/f_u(\lambda_j) \tag{27}$$

$$+ m^{-1}{\sum}'(j/p)^{2\theta}G_0(2\pi I_{\varepsilon j} - 1). \tag{28}$$

$\sup_{\Theta_2^1}|(26)|$, $\sup_{\Theta_2^1}|(27)| = o_p(1)$ follows from (19) and Lemmas 5.2(b) and 5.4. For (28), summation by parts gives

$$\begin{aligned}(28) &= G_0\left(\frac{m}{p}\right)^{2\theta}\frac{1}{m}\sum_{r=[\kappa m]}^{m-1}\left(\left(\frac{r}{m}\right)^{2\theta}-\left(\frac{r+1}{m}\right)^{2\theta}\right)\sum_{j=[\kappa m]}^{r}(2\pi I_{\varepsilon j}-1)\\ &\quad + G_0\left(\frac{m}{p}\right)^{2\theta}\frac{1}{m}\sum_{j=[\kappa m]}^{m}(2\pi I_{\varepsilon j}-1)\\ &= I(\theta)+II(\theta).\end{aligned}$$

As in [11], page 1637, write

$$2\pi I_{\varepsilon j}-1=\frac{1}{n}\sum_{t=1}^{n}(\varepsilon_t^2-1)+\frac{1}{n}\sum\sum_{s\neq t}\cos\{(s-t)\lambda_j\}\varepsilon_s\varepsilon_t,$$

from which it follows that

$$\begin{aligned}\sup_{\Theta_2^1}|I(\theta)| &\le \frac{C}{m}\sum_{r=[\kappa m]}^{m}\left|\sup_{\Theta_2^1}\left(\frac{r}{m}\right)^{2\theta}\right|\left|\frac{1}{n}\sum_{t=1}^{n}(\varepsilon_t^2-1)\right|\\ &\quad+\frac{C}{m}\sum_{r=[\kappa m]}^{m}\left|\sup_{\Theta_2^1}\left(\frac{r}{m}\right)^{2\theta}\right|\left|\frac{1}{rn}\right|\sum\sum_{s\neq t}\sum_{j=[\kappa m]}^{r}\cos\{(s-t)\lambda_j\}\varepsilon_s\varepsilon_t\Bigg|.\end{aligned}$$

From [11], (3.19) and (3.20), we have $n^{-1}\sum_1^n(\varepsilon_t^2-1)\xrightarrow{p}0$ and

$$E\left(\sum\sum_{s\neq t}\varepsilon_s\varepsilon_t\sum_{j=[\kappa m]}^{r}\cos\{(s-t)\lambda_j\}\right)^2=O(rn^2).$$

In conjunction with $\max_{[\kappa m]\le r\le m}\sup_{\Theta_2}(r/m)^{2\theta}=O(1)$, we obtain $\sup_{\Theta_2^1}|I(\theta)| = o_p(1)$. $\sup_{\Theta_2^1}|II(\theta)|=o_p(1)$ follows from a similar argument. Hence $\sup_{\Theta_2^1}|(28)| = o_p(1)$ and $\sup_{\Theta_2^1}|\Lambda_{2n}(\theta)|=o_p(1)$ follow, and we have established (23) for $\theta\in\Theta_2^1$.

For $\Theta_2^2=\{\theta: \frac{1}{2}\le\theta\le\frac{3}{2}\}$ define $Z_n(\theta)=\sum_{t=1}^n Y_t(\theta)\sim I(1-\theta)$ with $1-\theta\in[-\frac{1}{2},\frac{1}{2}]$. From Lemma 5.1(b) we have

$$w_{yj}=(1-e^{i\lambda_j})w_{zj}+(2\pi n)^{-1/2}e^{i\lambda_j}Z_n(\theta). \tag{29}$$

Define

$$\begin{aligned}D_{nj}(\theta)&=\lambda_j^{-\theta}(1-e^{i\lambda_j})D_n(e^{i\lambda_j};\theta-1),\\ \bar{U}_{nj}(\theta)&=\lambda_j^{-\theta}(1-e^{i\lambda_j})(2\pi n)^{-1/2}\widetilde{U}_{\lambda_j n}(\theta-1),\end{aligned}$$

and then applying (14) to $(Z_t(\theta),u_t)$ gives

$$\lambda_j^{-\theta}w_{yj}=D_{nj}(\theta)w_{uj}-\bar{U}_{nj}(\theta)+\lambda_j^{-\theta}(2\pi n)^{-1/2}e^{i\lambda_j}Z_n(\theta). \tag{30}$$

Since $\theta - 1 \ge -\frac{1}{2}$, from Lemma 5.2 we have, for $\theta \in \Theta_2^2$,

$$D_{nj}(\theta) = e^{-(\pi/2)\theta i} + O(\lambda_j) + O(j^{-1/2}) \qquad \text{uniformly in } \theta. \tag{31}$$

With a slight abuse of notation, rewrite

$$\begin{aligned} & m^{-1}\sum{}'(j/p)^{2\theta}(\lambda_j^{-2\theta} I_{yj} - G_0) \\ &\quad = m^{-1}\sum{}'(j/p)^{2\theta}[\lambda_j^{-2\theta} I_{yj} - |D_{nj}(\theta)|^2 I_{uj}] \\ &\qquad + m^{-1}\sum{}'(j/p)^{2\theta}[|D_{nj}(\theta)|^2 I_{uj} - G_0] \\ &\quad = \Lambda_{1n}(\theta) + \Lambda_{2n}(\theta). \end{aligned} \tag{32}$$

Therefore, (23) follows if, for $\theta \in \Theta_2^2$,

$$\Pr\Big(\inf_\theta \Lambda_{1n}(\theta) \le -2\delta G_0\Big) \to 0, \qquad \sup_\theta |\Lambda_{2n}(\theta)| = o_p(1) \qquad \text{as } n \to \infty. \tag{33}$$

$\sup_\theta |\Lambda_{2n}(\theta)| = o_p(1)$ follows straightforwardly from (31) and by the same argument as the one for $\theta \in \Theta_2^1$. For $\Lambda_{1n}(\theta)$, substituting (30) to the definition of $\Lambda_{1n}(\theta)$ gives

$$\Lambda_{1n}(\theta) = m^{-1}\sum{}'(j/p)^{2\theta}|\bar{U}_{nj}(\theta)|^2 \tag{34}$$

$$+ m^{-1}\sum{}'(j/p)^{2\theta}\lambda_j^{-2\theta}(2\pi n)^{-1} Z_n^2 \tag{35}$$

$$- 2\,\mathrm{Re}\Big[m^{-1}\sum{}'(j/p)^{2\theta} D_{nj}(\theta)^* w_{uj}^* \bar{U}_{nj}(\theta)\Big] \tag{36}$$

$$- 2\,\mathrm{Re}\Big[m^{-1}\sum{}'(j/p)^{2\theta} \bar{U}_{nj}(\theta)\lambda_j^{-\theta}(2\pi n)^{-1/2} e^{i\lambda_j} Z_n(\theta)\Big] \tag{37}$$

$$+ 2\,\mathrm{Re}\Big[m^{-1}\sum{}'(j/p)^{2\theta} D_{nj}(\theta)^* w_{uj}^* \lambda_j^{-\theta}(2\pi n)^{-1/2} e^{i\lambda_j} Z_n(\theta)\Big]. \tag{38}$$

Equation (34) is almost surely nonnegative. Lemma 5.3 gives

$$E\sup_\theta |\bar{U}_{nj}(\theta)|^2 = O(j^{-1}(\log n)^2), \tag{39}$$

and hence $\sup_\theta |(36)| = o_p(1)$ follows from (39) and Lemma 5.4. Therefore, (33) and hence (23) follow if, for any $\zeta > 0$,

$$\Pr\Big(\inf_\theta [(35) + (37) + (38)] \le -\zeta\Big) \to 0 \qquad \text{as } n \to \infty. \tag{40}$$

We proceed to show (40). First, there exists $\eta > 0$ such that, uniformly in $\theta$,

$$(35) = p^{-2\theta}(2\pi)^{-2\theta-1} n^{2\theta-1} Z_n(\theta)^2 m^{-1} \sum{}' 1 \geq \eta (m^{-\theta} n^{\theta-1/2} Z_n(\theta))^2.$$

From (39) and Lemma 5.4, we have, uniformly in $\theta$,

$$(37) = m^{-\theta} n^{\theta-1/2} Z_n(\theta) \cdot O_p(m^{-1/2} \log n).$$

For (38), it follows from (31), $e^{i\lambda_j} = 1 + O(\lambda_j)$ and Lemmas 5.4 and 5.6 that

$$\begin{aligned} & m^{-1} \sum{}' (j/p)^{2\theta} D_{nj}(\theta)^* w_{uj}^* \lambda_j^{-\theta} (2\pi n)^{-1/2} e^{i\lambda_j} Z_n(\theta) \\ & \quad = (2\pi n)^{-1/2} Z_n(\theta) e^{(\pi/2)\theta i} m^{-1} \sum{}' (j/p)^{2\theta} w_{uj}^* \lambda_j^{-\theta} \\ & \qquad + (2\pi n)^{-1/2} Z_n(\theta) m^{-1} \sum{}' (j/p)^{2\theta} w_{uj}^* \lambda_j^{-\theta} [O(\lambda_j) + O(j^{-1/2})] \\ & \quad = m^{-\theta} n^{\theta-1/2} Z_n(\theta)[O_p(m^{-1/2} \log m) + O_p(mn^{-1})]. \end{aligned} \tag{41}$$

Therefore, we can write

$$(37) + (38) = m^{-\theta} n^{\theta-1/2} Z_n(\theta) \cdot R_n(\theta, \omega), \tag{42}$$

where $\omega$ denotes an element of the sample space $\Omega$, and

$$\sup_\theta |R_n(\theta, \omega)| = O_p(k_n), \qquad k_n = m^{-1/2} \log n + mn^{-1} \to 0. \tag{43}$$

Before showing (40), define

$$\begin{aligned} \Omega_1 &= \{(\omega, \theta) \in \Omega \times \Theta : m^{-\theta} n^{\theta-1/2} |Z_n(\theta)| < k_n \log m\}, \\ \Omega_2 &= \{(\omega, \theta) \in \Omega \times \Theta : m^{-\theta} n^{\theta-1/2} |Z_n(\theta)| \geq k_n \log m\}, \end{aligned}$$

where $\Theta$ is the domain of $\theta$ ($\Theta_2^1$ in this case), so that $\Omega_1 \cup \Omega_2 = \Omega \times \Theta$. Then

$$\begin{aligned} & \{(\omega, \theta) : \eta (m^{-\theta} n^{\theta-1/2} Z_n(\theta))^2 - |m^{-\theta} n^{\theta-1/2} Z_n(\theta) \cdot R_n(\theta, \omega)| \leq -\zeta\} \\ & \quad = \{(\omega, \theta) : (\eta (m^{-\theta} n^{\theta-1/2} Z_n(\theta))^2 \\ & \qquad\qquad - |m^{-\theta} n^{\theta-1/2} Z_n(\theta) \cdot R_n(\theta, \omega)| \leq -\zeta) \cap \Omega_1\} \\ & \qquad \cup \{(\omega, \theta) : (\eta (m^{-\theta} n^{\theta-1/2} Z_n(\theta))^2 \\ & \qquad\qquad - |m^{-\theta} n^{\theta-1/2} Z_n(\theta) \cdot R_n(\theta, \omega)| \leq -\zeta) \cap \Omega_2\} \\ & \quad \subseteq \{(\omega, \theta) : \eta (m^{-\theta} n^{\theta-1/2} Z_n(\theta))^2 - k_n \log m |R_n(\theta, \omega)| \leq -\zeta\} \\ & \qquad \cup \{(\omega, \theta) : m^{-\theta} n^{\theta-1/2} |Z_n(\theta)| [\eta k_n \log m - |R_n(\theta, \omega)|] \leq -\zeta\} \\ & \quad \subseteq \{(\omega, \theta) : k_n \log m |R_n(\theta, \omega)| \geq \zeta\} \cup \{(\omega, \theta) : \eta k_n \log m - |R_n(\theta, \omega)| \leq 0\}. \end{aligned}$$

Therefore,

$$
\begin{aligned}
&\left\{\omega : \inf_\theta [\eta(m^{-\theta} n^{\theta-1/2} Z_n(\theta))^2 - |m^{-\theta} n^{\theta-1/2} Z_n(\theta) \cdot R_n(\theta,\omega)|] \le -\zeta\right\} \\
&\quad \subseteq \left\{\omega : \sup_\theta k_n \log m |R_n(\theta,\omega)| \ge \delta G_0\right\} \\
&\qquad \cup \left\{\omega : \eta k_n \log m - \sup_\theta |R_n(\theta,\omega)| \le 0\right\},
\end{aligned}
$$

and it follows that

$$
\begin{aligned}
&\Pr\left(\inf_\theta [\eta(m^{-\theta} n^{\theta-1/2} Z_n(\theta))^2 - |m^{-\theta} n^{\theta-1/2} Z_n \cdot R_n(\theta,\omega)|] \le -\zeta\right) \\
&\quad \le \Pr\left(k_n \log m \sup_\theta |R_n(\theta,\omega)| \ge \zeta\right) \\
&\qquad + \Pr\left(\eta k_n \log m - \sup_\theta |R_n(\theta,\omega)| \le 0\right) \\
&\quad \to 0,
\end{aligned}
$$

because $\sup_\theta |R_n(\theta,\omega)| = O_p(k_n)$, and $k_n^2 \log m \to 0$ from Assumption 4. Therefore (40) follows, and hence (23) holds for $\theta \in \Theta_2^2$.

For $\Theta_2^3 = \{\theta : -\frac{3}{2} \le \theta \le -\frac{1}{2}\}$, from Lemma 5.1 we have

$$
w_{yj} = (1 - e^{i\lambda_j})^{-1} w_{\Delta yj} - (1 - e^{i\lambda_j})^{-1} (2\pi n)^{-1/2} e^{i\lambda_j} Y_n(\theta), \tag{44}
$$

where $\Delta Y_t(\theta) \sim I(-\theta - 1)$. With a slight abuse of notation, define

$$
\begin{aligned}
D_{nj}(\theta) &= \lambda_j^{-\theta} (1 - e^{i\lambda_j})^{-1} D_n(e^{i\lambda_j}; \theta + 1), \\
\bar{U}_{nj}(\theta) &= \lambda_j^{-\theta} (1 - e^{i\lambda_j})^{-1} (2\pi n)^{-1/2} \widetilde{U}_{\lambda_j n}(\theta + 1).
\end{aligned}
$$

Then, applying (14) to $(\Delta Y_t(\theta), u_t)$ gives

$$
\lambda_j^{-\theta} w_{yj} = D_{nj}(\theta) w_{uj} - \bar{U}_{nj}(\theta) + \lambda_j^{-\theta} (2\pi n)^{-1/2} e^{i\lambda_j} (1 - e^{i\lambda_j})^{-1} Y_n(\theta). \tag{45}
$$

$D_{nj}(\theta)$ and $\bar{U}_{nj}(\theta)$ satisfy (31) and (39) for $\theta \in \Theta_2^3$, because $-\theta - 1 \in [-\frac{1}{2}, \frac{1}{2}]$. Using the decomposition (32) and the same argument as the one for $\theta \in \Theta_2^2$, we obtain

$$
\begin{aligned}
&m^{-1} {\sum}' (j/p)^{2\theta} (\lambda_j^{-2\theta} I_{yj} - G_0) \\
&\quad = m^{-1} {\sum}' (j/p)^{2\theta} [\lambda_j^{-2\theta} I_{yj} - |D_{nj}(\theta)|^2 I_{uj}] + o_p(1),
\end{aligned}
$$

where $o_p(1)$ is uniform in $\theta \in \Theta_2^3$. Using (45), rewrite the first term on the right-hand side as

$$m^{-1}\sum{}'(j/p)^{2\theta}|\bar{U}_{nj}(\theta)|^2 \tag{46}$$

$$+\,m^{-1}\sum{}'(j/p)^{2\theta}\lambda_j^{-2\theta}(2\pi n)^{-1}|1-e^{i\lambda_j}|^{-2}Y_n(\theta)^2 \tag{47}$$

$$-\,2\,\mathrm{Re}\left[m^{-1}\sum{}'(j/p)^{2\theta}D_{nj}(\theta)^* w_{uj}^*\bar{U}_{nj}(\theta)\right] \tag{48}$$

$$-\,2\,\mathrm{Re}\left[m^{-1}\sum{}'(j/p)^{2\theta}\bar{U}_{nj}(\theta)\lambda_j^{-\theta}(2\pi n)^{-1/2}e^{i\lambda_j}(1-e^{i\lambda_j})^{-1}Y_n(\theta)\right] \tag{49}$$

$$\begin{aligned}&+\,2\,\mathrm{Re}\Bigg[m^{-1}\sum{}'(j/p)^{2\theta}D_{nj}(\theta)^*\\&\qquad\times w_{uj}^*\lambda_j^{-\theta}(2\pi n)^{-1/2}e^{i\lambda_j}(1-e^{i\lambda_j})^{-1}Y_n(\theta)\Bigg].\end{aligned} \tag{50}$$

Equation (46) is almost surely nonnegative. Because $D_{nj}(\theta)$ and $\bar{U}_{nj}(\theta)$ satisfy (31) and (39), it follows from a decomposition similar to (41) and Lemmas 5.4 and 5.6 that $\sup_\theta|(48)| = o_p(1)$ and $(49)+(50) = m^{-\theta-1}n^{\theta+1/2}Y_n(\theta)\times O_p(m^{-1/2}\log n + mn^{-1})$. Finally, (47) is equal to

$$\begin{aligned}&p^{-2\theta}n^{2\theta-1}(2\pi)^{-2\theta-1}m^{-1}\sum{}'|1-e^{i\lambda_j}|^{-2}Y_n(\theta)^2\\&\quad = p^{-2\theta}n^{2\theta-1}(2\pi)^{-2\theta-1}Y_n(\theta)^2 m^{-1}\sum{}'\lambda_j^{-2}(1+o(1))\\&\quad \ge \eta m^{-2\theta-2}n^{2\theta+1}Y_n(\theta)^2,\end{aligned} \tag{51}$$

for some $\eta > 0$. Therefore we can apply the argument following (42) with slight changes to show (23) for $\theta \in \Theta_2^3$.

For $\Theta_2^4 = \{\theta : \frac{3}{2} \le \theta \le \frac{5}{2}\}$, by applying (29) twice and (14), we obtain

$$\begin{aligned}\lambda_j^{-\theta}w_{yj} &= D_{nj}(\theta)w_{uj} - \bar{U}_{nj}(\theta)\\&\qquad + \lambda_j^{-\theta}(2\pi n)^{-1/2}e^{i\lambda_j}\left[(1-e^{i\lambda_j})\sum_{t=1}^{n}Z_t(\theta) + Z_n(\theta)\right],\end{aligned}$$

where

$$D_{nj}(\theta) = \lambda_j^{-\theta}(1-e^{i\lambda_j})^2 D_n(e^{i\lambda_j};\theta-2),$$

$$\bar{U}_{nj}(\theta) = \lambda_j^{-\theta}(1-e^{i\lambda_j})^2(2\pi n)^{-1/2}\widetilde{U}_{\lambda_j n}(\theta-2),$$

and $D_{nj}(\theta)$ and $\bar{U}_{nj}(\theta)$ satisfy (31) and (39). We proceed to evaluate the terms in $m^{-1}\sum'(j/p)^{2\theta}\lambda_j^{-2\theta}I_{yj}$. First, observe that

$$
\begin{aligned}
&m^{-1}\sum{}'(j/p)^{2\theta}\lambda_j^{-2\theta}(2\pi n)^{-1}\left|(1-e^{i\lambda_j})\sum_{t=1}^{n}Z_t(\theta)+Z_n(\theta)\right|^2 \\
&\qquad = p^{-2\theta}n^{2\theta-1}(2\pi)^{-2\theta-1}m^{-1}\sum{}'\left|(1-e^{i\lambda_j})\sum_{t=1}^{n}Z_t(\theta)+Z_n(\theta)\right|^2 .
\end{aligned}\tag{52}
$$

By applying Lemma 5.10(a) with $Q_3 = Q_2 = 0$, $Q_1 = \sum_1^n Z_t(\theta)$ and $Q_0 = Z_n(\theta)$, there exists $\eta > 0$ such that, for sufficiently large $n$,

$$
(52) \geq \eta m^{-2\theta+2}n^{2\theta-3}\left(\sum_{t=1}^{n}Z_t(\theta)\right)^2 + \eta m^{-2\theta}n^{2\theta-1}Z_n(\theta)^2 = \Lambda_{3n}(\theta)
$$

uniformly in $\theta$. Of the other terms in $m^{-1}\sum'(j/p)^{2\theta}\lambda_j^{-2\theta}I_{yj}$, the terms involving the cross products of $w_{uj}, \bar{U}_{nj}(\theta)$ and $(1-e^{i\lambda_j})\sum_1^n Z_t(\theta)+Z_n(\theta)$ are dominated by $\Lambda_{3n}(\theta)$. For instance, proceeding as in (41) gives

$$
\begin{aligned}
&m^{-1}\sum{}'(j/p)^{2\theta}D_{nj}(\theta)w_{uj}\lambda_j^{-\theta}(2\pi n)^{-1/2}e^{-i\lambda_j}\left[(1-e^{-i\lambda_j})\sum_{t=1}^{n}Z_t(\theta)+Z_n(\theta)\right] \\
&\quad = m^{-\theta+1}n^{\theta-3/2}\sum_{t=1}^{n}Z_t(\theta)\cdot O_p(m^{-1/2}\log n + n^{-1}m) \\
&\qquad + m^{-\theta}n^{\theta-1/2}Z_n(\theta)\cdot O_p(m^{-1/2}\log n + n^{-1}m),
\end{aligned}
$$

where the $O_p(\cdot)$ terms are uniform in $\theta$. Therefore, the terms in $m^{-1}\sum'(j/p)^{2\theta}\times[\lambda_j^{-2\theta}I_{yj} - |D_{nj}(\theta)|^2 I_{uj}]$ are either $o_p(1)$ or nonnegative or dominated by $\Lambda_{3n}(\theta)$. We obtain $\sup_\theta|m^{-1}\sum'(j/p)^{2\theta}[|D_{nj}(\theta)|^2 I_{uj} - G_0]| = o_p(1)$ by using (31) and proceeding as in (26)–(28) and the following argument, and thus (23) follows for $\theta \in \Theta_2^4$.

Since $|\theta| \leq \Delta_2 - \Delta_1 \leq \frac{9}{2}$, the proof is completed by showing (23) for the remaining subsets of $\Theta_2$:

$$
\begin{aligned}
\Theta_2^5 &= \{\theta : -\tfrac{5}{2} \leq \theta \leq -\tfrac{3}{2}\},\\
\Theta_2^6 &= \{\theta : \tfrac{7}{2} \leq \theta \leq \tfrac{5}{2}\},\\
\Theta_2^7 &= \{\theta : -\tfrac{7}{2} \leq \theta \leq -\tfrac{5}{2}\},\\
\Theta_2^8 &= \{\theta : \tfrac{9}{2} \leq \theta \leq \tfrac{7}{2}\},\\
\Theta_2^9 &= \{\theta : -\tfrac{9}{2} \leq \theta \leq -\tfrac{7}{2}\}.
\end{aligned}
$$

Applying (29) or (44) repeatedly and (14) gives the required result for $\Theta_2^{\cdot}$. For instance, for $\Theta_2^9 = \{\theta : -\frac{9}{2} \le \theta \le -\frac{7}{2}\}$, applying (44) four times and then (14), we have

$$\lambda_j^{-\theta} w_{yj} = D_{nj}(\theta) w_{uj} - \bar{U}_{nj}(\theta) - \lambda_j^{-\theta}(2\pi n)^{-1/2} e^{i\lambda_j} W_{nj},$$

where

$$\begin{aligned} D_{nj}(\theta) &= \lambda_j^{-\theta}(1 - e^{i\lambda_j})^{-4} D_n(e^{i\lambda_j}; \theta + 4), \\ \bar{U}_{nj}(\theta) &= \lambda_j^{-\theta}(1 - e^{i\lambda_j})^{-4}(2\pi n)^{-1/2} \widetilde{U}_{\lambda_j n}(\theta + 4), \\ W_{nj} &= (1 - e^{i\lambda_j})^{-4} \Delta^3 Y_n(\theta) - (1 - e^{i\lambda_j})^{-3} \Delta^2 Y_n(\theta) \\ &\quad - (1 - e^{i\lambda_j})^{-2} \Delta Y_n(\theta) - (1 - e^{i\lambda_j})^{-1} Y_n(\theta), \end{aligned}$$

and $D_{nj}(\theta)$ and $\bar{U}_{nj}(\theta)$ satisfy (31) and (39). We can easily obtain

$$\begin{aligned} & m^{-1} \sum{}' (j/p)^{2\theta} (\lambda_j^{-2\theta} I_{yj} - G_0) \\ & \qquad = m^{-1} \sum{}' (j/p)^{2\theta} [\lambda_j^{-2\theta} I_{yj} - |D_{nj}(\theta)|^2 I_{uj}] + o_p(1), \end{aligned}$$

where $o_p(1)$ is uniform in $\theta \in \Theta_2^9$. For the first term on the right-hand side, from Lemma 5.10(b) we have, for large $n$ and $\eta > 0$,

$$\begin{aligned} & m^{-1} \sum{}' (j/p)^{2\theta} \lambda_j^{-2\theta} (2\pi n)^{-1} |W_{nj}|^2 \\ & \quad = (2\pi)^{-2\theta-1} n^{2\theta-1} p^{-2\theta} m^{-1} \sum{}' |W_{nj}|^2 \\ & \quad \ge \eta n^{2\theta-1} m^{-2\theta} \begin{bmatrix} m^{-8} n^8 (\Delta^3 Y_n(\theta))^2 + m^{-6} n^6 (\Delta^2 Y_n(\theta))^2 \\ + m^{-4} n^4 (\Delta Y_n(\theta))^2 + m^{-2} n^2 Y_n(\theta)^2 \end{bmatrix}, \end{aligned} \tag{53}$$

uniformly in $\theta$. The terms involving the cross products between $w_{uj}$, $\bar{U}_{nj}(\theta)$ and $W_{nj}$ are dominated by (53). The other terms in $m^{-1} \sum' (j/p)^{2\theta} [\lambda_j^{-2\theta} I_{yj} - |D_{nj}(\theta)|^2 I_{uj}]$ are either $o_p(1)$ or almost surely nonnegative, and hence (23) follows.

4.2. *Proof of asymptotic normality.* Theorem 2.1 holds under the current conditions and implies that with probability approaching 1, as $n \to \infty$ $\widehat{d}$ satisfies

$$0 = R'(\widehat{d}) = R'(d_0) + R''(\bar{d})(\widehat{d} - d_0), \tag{54}$$

where $|\bar{d} - d_0| \leq |\widehat{d} - d_0|$. From the fact that

$$\frac{\partial}{\partial d} w_{\Delta^d xj} = \frac{\partial}{\partial d} \frac{1}{\sqrt{2\pi n}} \sum_{t=1}^{n} e^{i\lambda_j t}(1-L)^d X_t$$

$$= \frac{1}{\sqrt{2\pi n}} \sum_{t=1}^{n} e^{i\lambda_j t} \log(1-L)(1-L)^d X_t,$$

$$\frac{\partial^2}{\partial d^2} w_{\Delta^d xj} = \frac{1}{\sqrt{2\pi n}} \sum_{t=1}^{n} e^{i\lambda_j t} (\log(1-L))^2 (1-L)^d X_t,$$

we obtain

$$R''(d) = \frac{\widehat{G}_2(d)\widehat{G}(d) - \widehat{G}_1^2(d)}{\widehat{G}^2(d)} = \frac{\widetilde{G}_2(d)\widetilde{G}_0(d) - \widetilde{G}_1^2(d)}{\widetilde{G}_0^2(d)},$$

where

$$\widehat{G}_1(d) = \frac{1}{m} \sum_{j=1}^{m} \frac{\partial}{\partial d}[w_{\Delta^d xj} w^*_{\Delta^d xj}] = \frac{1}{m} \sum_{j=1}^{m} 2\,\mathrm{Re}[w_{\log(1-L)\Delta^d xj} w^*_{\Delta^d xj}],$$

$$\widehat{G}_2(d) = \frac{1}{m} \sum_{j=1}^{m} \frac{\partial^2}{\partial d^2}[w_{\Delta^d xj} w^*_{\Delta^d xj}] = \frac{1}{m} \sum_{j=1}^{m} W_x(L, d, j),$$

$$W_x(L, d, j) = 2\,\mathrm{Re}[w_{(\log(1-L))^2 \Delta^d xj} w^*_{\Delta^d xj}] + 2I_{\log(1-L)\Delta^d xj},$$

$$\widetilde{G}_0(d) = \frac{1}{m} \sum_{j=1}^{m} j^{2\theta} \lambda_j^{-2\theta} I_{yj},$$

$$\widetilde{G}_1(d) = \frac{1}{m} \sum_{j=1}^{m} j^{2\theta} \lambda_j^{-2\theta} 2\,\mathrm{Re}[w_{\log(1-L)yj} w^*_{yj}],$$

$$\widetilde{G}_2(d) = \frac{1}{m} \sum_{j=1}^{m} j^{2\theta} \lambda_j^{-2\theta} W_y(L, 0, j),$$

and $\theta = d - d_0$ and $Y_t(\theta) = (1-L)^d X_t = (1-L)^\theta u_t I\{t \geq 1\}$, as defined in the proof of Theorem 2.1. Fix $\varepsilon > 0$ and let $M = \{d : (\log n)^4 |d - d_0| < \varepsilon\}$. From (9) in the proof of Theorem 2.1 we have

$$\Pr(\bar{d} \notin M) \leq \left(\inf_{\Theta_1 \setminus M} S(d) \leq 0\right) + o(1).$$

Hence, in view of (10), $\Pr(\bar{d} \notin M)$ tends to zero if

$$\sup_{\Theta_1} |A(d)/B(d)| = o_p((\log n)^{-8}), \tag{55}$$

where $A(d)$ and $B(d)$ are defined in (12) in the proof of Theorem 2.1. From Assumption 1′, (18) is strengthened to

$$|\lambda_j^{-2\theta}|D_n(e^{i\lambda_j};\theta)|^2 - G_0/f_u(\lambda_j)| = O(\lambda_j^{\beta}) + O(j^{-1/2}), \qquad j = 1,\dots,m. \tag{56}$$

Therefore, proceeding as in the proof of Theorem 2.1, we obtain

$$\sum_{r=1}^{m}\left(\frac{r}{m}\right)^{2\Delta}\frac{1}{r^2}\sup_{\Theta_1}\left|\sum_{j=1}^{r}[\lambda_j^{-2\theta}I_{yj} - 2\pi G_0 I_{\varepsilon j}]\right| = O_p(m^{\beta}n^{-\beta} + m^{-\Delta}(\log n)^2).$$

Robinson ([11], (4.9), page 1643) shows

$$\sum_{j=1}^{r}(2\pi I_{\varepsilon j} - 1) = O_p(r^{1/2}) \qquad \text{as } n \to \infty \text{ for } 1 \le r \le m, \tag{57}$$

and it follows that $\sum_1^m (r/m)^{2\Delta} r^{-2}|\sum_1^r (2\pi I_{\varepsilon j} - 1)| = O(m^{-2\Delta}\log m)$. Applying the same argument to the second term of (16), we obtain $\sup_{\Theta_1}|A(d)| = o_p((\log n)^{-8})$, and (55) follows in view of (11). Thus we assume $d \in M$ in the following discussion of $\widetilde{G}_k(d)$.

Now we derive the approximation of $\widetilde{G}_k(d)$ for $k = 0,1,2$. For $\widetilde{G}_0(d)$ observe that

$$\begin{aligned} E\sup_{\theta\in M}&|\lambda_j^{-2\theta}I_{yj} - I_{uj}| \\ &\le E\sup_{\theta\in M}|\lambda_j^{-2\theta}I_{yj} - \lambda_j^{-2\theta}|D_n(e^{i\lambda_j};\theta)|^2 I_{uj}| \\ &\quad + E\sup_{\theta\in M}|\lambda_j^{-2\theta}|D_n(e^{i\lambda_j};\theta)|^2 - 1|I_{uj} \\ &= O(j^{-1/2}(\log n)^2 + j^2 n^{-2}), \qquad j = 1,\dots,m, \end{aligned} \tag{58}$$

where the third line follows from (21) and Lemma 5.2. Since $|j^{2\theta} - 1|/|2\theta| \le (\log j)n^{2|\theta|} \le (\log j)n^{1/\log n} = e\log j$ on $M$, we have

$$\begin{aligned} &\sup_M |j^{2\theta} - 1| = O((\log n)^{-3}), \\ &\sup_M |j^{2\theta}| = O(1), \qquad j = 1,\dots,m. \end{aligned} \tag{59}$$

Therefore, in view of (58) and $EI_{uj} = O(1)$ [following from (19)], we obtain

$$\begin{aligned} \sup_M&\left|\widetilde{G}_0(d) - \frac{1}{m}\sum_{j=1}^{m}I_{uj}\right| \\ &\le \sup_M\left|\frac{1}{m}\sum_{j=1}^{m}j^{2\theta}[\lambda_j^{-2\theta}I_{yj} - I_{uj}]\right| + \sup_M\left|\frac{1}{m}\sum_{j=1}^{m}(j^{2\theta} - 1)I_{uj}\right| \\ &= o_p((\log n)^{-2}). \end{aligned}$$

For $\widetilde{G}_1(d)$, from (14) and Lemma 5.9 we have

$$\begin{aligned}
&\lambda_j^{-2\theta} w_{\log(1-L)yj} w_{yj}^* + J_n(e^{i\lambda_j}) I_{uj} \\
&\quad = J_n(e^{i\lambda_j})[1 - \lambda_j^{-2\theta} |D_n(e^{i\lambda_j};\theta)|^2] I_{uj} \\
&\qquad - J_n(e^{i\lambda_j}) \lambda_j^{-\theta} D_n(e^{i\lambda_j};\theta) w_{uj} \cdot \lambda_j^{-\theta} (2\pi n)^{-1/2} \widetilde{U}_{\lambda_j n}(\theta)^* \\
&\qquad - \lambda_j^{-\theta} D_n(e^{i\lambda_j};\theta)^* w_{uj}^* \cdot \lambda_j^{-\theta} (2\pi n)^{-1/2} V_{nj}(\theta) \\
&\qquad - \lambda_j^{-2\theta} (2\pi n)^{-1} \widetilde{U}_{\lambda_j n}(\theta)^* V_{nj}(\theta).
\end{aligned}$$

Then, since $J_n(e^{i\lambda_j}) = O(\log n)$, it follows from (59) and Lemmas 5.2, 5.3 and 5.9 that

$$\frac{1}{m} \sum_{j=1}^{m} \sup_M j^{2\theta} |\operatorname{Re}[\lambda_j^{-2\theta} w_{\log(1-L)yj} w_{yj}^* + J_n(e^{i\lambda_j}) I_{uj}]| = o_p((\log n)^{-1}).$$

In conjunction with (59), $J_n(e^{i\lambda_j}) = O(\log n)$ and $EI_{uj} = O(1)$, it follows that

$$\begin{aligned}
&\sup_M \left| \widetilde{G}_1(d) + \frac{1}{m} \sum_{j=1}^{m} 2\operatorname{Re}[J_n(e^{i\lambda_j})] I_{uj} \right| \\
&\quad = \sup_M \left| \frac{1}{m} \sum_{j=1}^{m} (1 - j^{2\theta}) 2\operatorname{Re}[J_n(e^{i\lambda_j})] I_{uj} \right| + o_p((\log n)^{-1}) \\
&\quad = o_p((\log n)^{-1}).
\end{aligned}$$

For $\widetilde{G}_2(d)$, the same line of argument as above with Lemma 5.9(c) gives

$$\begin{aligned}
&\sup_M \left| \widetilde{G}_2(d) - \frac{1}{m} \sum_{j=1}^{m} \{2\operatorname{Re}[J_n(e^{i\lambda_j})^2] + 2J_n(e^{i\lambda_j}) J_n(e^{i\lambda_j})^*\} I_{uj} \right| \\
&\quad = \sup_M \left| \widetilde{G}_2(d) - \frac{1}{m} \sum_{j=1}^{m} 4\{\operatorname{Re}[J_n(e^{i\lambda_j})]\}^2 I_{uj} \right| \\
&\quad = o_p(1).
\end{aligned}$$

From (19) and Assumption $1'$, we obtain

$$\begin{aligned}
E|I_{uj} - G_0 I_{\varepsilon j}| &\le E|I_{uj} - |C(e^{i\lambda_j})|^2 I_{\varepsilon j}| + E2\pi |f_u(\lambda_j) - f_u(0)| I_{\varepsilon j} \\
&= O(j^{-1/2}(\log(j+1)) + j^{\beta} n^{-\beta}), \qquad j = 1, \dots, m.
\end{aligned}$$

Therefore, in view of $J_n(e^{i\lambda_j}) = O(\log n)$, $EI_{\varepsilon j} = 1$, and $Cov(I_{\varepsilon j}, I_{\varepsilon k}) = O(1)$ if $j = k$ and $O(n^{-1})$ if $j \ne k$, we have

$$\widetilde{G}_0(\bar{d}) = m^{-1} \sum_{j=1}^{m} I_{uj} + o_p((\log n)^{-2})$$

$$
\begin{aligned}
&= G_0 m^{-1} \sum_{j=1}^{m} I_{\varepsilon j} + o_p((\log n)^{-2}) \\
&= G_0 + o_p((\log n)^{-2}), \\
\widetilde{G}_1(\bar{d}) &= -2m^{-1} \sum_{j=1}^{m} \mathrm{Re}[J_n(e^{i\lambda_j})] I_{uj} + o_p((\log n)^{-1}) \\
&= -G_0 m^{-1} \sum_{j=1}^{m} 2\,\mathrm{Re}[J_n(e^{i\lambda_j})] I_{\varepsilon j} + o_p((\log n)^{-1}) \\
&= -G_0 m^{-1} \sum_{j=1}^{m} 2\,\mathrm{Re}[J_n(e^{i\lambda_j})] + o_p((\log n)^{-1})
\end{aligned}
$$

and

$$
\begin{aligned}
\widetilde{G}_2(\bar{d}) &= m^{-1} \sum_{j=1}^{m} 4\{\mathrm{Re}[J_n(e^{i\lambda_j})]\}^2 I_{uj} + o_p(1) \\
&= G_0 m^{-1} \sum_{j=1}^{m} 4\{\mathrm{Re}[J_n(e^{i\lambda_j})]\}^2 I_{\varepsilon j} + o_p(1) \\
&= G_0 m^{-1} \sum_{j=1}^{m} 4\{\mathrm{Re}[J_n(e^{i\lambda_j})]\}^2 + o_p(1).
\end{aligned}
$$

It follows that

$$
\begin{aligned}
R''(\bar{d}) &= [\widetilde{G}_2(\bar{d})\widetilde{G}_0(\bar{d}) - \widetilde{G}_1^2(\bar{d})][\widetilde{G}_0(\bar{d})]^{-2} \\
&= \frac{G_0^2 m^{-1} \sum_1^m 4\{\mathrm{Re}[J_n(e^{i\lambda_j})]\}^2 - \{G_0 m^{-1} \sum_1^m 2\,\mathrm{Re}[J_n(e^{i\lambda_j})]\}^2 + o_p(1)}{\{G_0 + o_p((\log n)^{-2})\}^2} \qquad (60) \\
&= 4m^{-1} \sum_{j=1}^{m} \{\mathrm{Re}[J_n(e^{i\lambda_j})]\}^2 - 4\left\{ m^{-1} \sum_{j=1}^{m} \mathrm{Re}[J_n(e^{i\lambda_j})] \right\}^2 + o_p(1).
\end{aligned}
$$

From Lemma 5.8(a) and a routine calculation, we obtain

$$
\begin{aligned}
m^{-1} \sum_{j=1}^{m} \{\mathrm{Re}[J_n(e^{i\lambda_j})]\}^2 &= m^{-1} \sum_{j=1}^{m} (\log \lambda_j)^2 + o(1), \\
\left\{ m^{-1} \sum_{j=1}^{m} \mathrm{Re}[J_n(e^{i\lambda_j})] \right\}^2 &= \left( m^{-1} \sum_{j=1}^{m} \log \lambda_j \right)^2 + o(1).
\end{aligned}
$$

Therefore, $\frac{1}{4}$ times (60) is, apart from $o_p(1)$ terms,

$$\begin{aligned} m^{-1}\sum_{j=1}^{m}(\log\lambda_j)^2 &- \left(m^{-1}\sum_{j=1}^{m}\log\lambda_j\right)^2 \\ &= m^{-1}\sum_{j=1}^{m}(\log j)^2 - \left(m^{-1}\sum_{j=1}^{m}\log j\right)^2 \to 1, \end{aligned}$$

and $R''(\bar{d}) = 4 + o_p(1)$ follows.

Now we find the limit distribution of $m^{1/2}R'(d_0)$. In view of Lemma 5.9(b), $E|w_{uj} - C(e^{i\lambda_j})w_{\varepsilon j}|^2 = O(j^{-1}\log(j+1))$ [see (19)] and $E|\widetilde{J}_{n\lambda_j}(e^{i\lambda_j}L)\varepsilon_n|^2 = O(nj^{-1})$ [see (77)], we obtain

$$\begin{aligned} -w_{\log(1-L)uj}w_{uj}^* &\\ &= [J_n(e^{i\lambda_j})w_{uj} + r_{nj}]w_{uj}^* \\ &\quad - C(1)(2\pi n)^{-1/2}\widetilde{J}_{n\lambda_j}(e^{-i\lambda_j}L)\varepsilon_n C(e^{i\lambda_j})^* w_{\varepsilon j}^* \\ &\quad - C(1)(2\pi n)^{-1/2}\widetilde{J}_{n\lambda_j}(e^{-i\lambda_j}L)\varepsilon_n[w_{uj}^* - C(e^{i\lambda_j})^* w_{\varepsilon j}^*] \\ &= J_n(e^{i\lambda_j})I_{uj} - C(1)(2\pi n)^{-1/2}\widetilde{J}_{n\lambda_j}(e^{-i\lambda_j}L)\varepsilon_n C(e^{i\lambda_j})^* w_{\varepsilon j}^* + R_{nj}, \end{aligned}$$

where $r_{nj}$ is defined in Lemma 5.9(b), and $E|j^{1/2}R_{nj}| = o(1) + O(j^{-1/2}\log m)$ as $n \to \infty$. It follows that $m^{1/2}\widehat{G}_1(d_0)$ is equal to

$$-m^{-1/2}\sum_{j=1}^{m} 2\,\mathrm{Re}[J_n(e^{i\lambda_j})]I_{uj} \tag{61}$$

$$+ C(1)m^{-1/2}\sum_{j=1}^{m} 2\,\mathrm{Re}[(2\pi n)^{-1/2}\widetilde{J}_{n\lambda_j}(e^{-i\lambda_j}L)\varepsilon_n C(e^{i\lambda_j})^* w_{\varepsilon j}^*] \tag{62}$$

$$+ o_p(1) + O_p(m^{-1/2}(\log m)^2).$$

From Lemma 5.8(a) we have

$$(61) = 2m^{-1/2}\sum_{j=1}^{m}(\log\lambda_j)I_{uj} + O_p(m^{5/2}n^{-2}) + O_p(m^{-1/2}\log m).$$

For (62), in view of the fact that

$$w_{\varepsilon j}^* = (2\pi n)^{-1/2}\sum_{p=1}^{n} e^{-ip\lambda_j}\varepsilon_p = (2\pi n)^{-1/2}\sum_{q=0}^{n-1} e^{iq\lambda_j}\varepsilon_{n-q},$$

we obtain the decomposition

$$
\begin{aligned}
&m^{-1/2}\sum_{j=1}^{m}(2\pi n)^{-1/2}\widetilde{J}_{n\lambda_j}(e^{-i\lambda_j}L)\varepsilon_n C(e^{i\lambda_j})^* w_{\varepsilon j}^* \\
&\quad = m^{-1/2}\sum_{j=1}^{m} C(e^{i\lambda_j})^*(2\pi n)^{-1}\left(\sum_{p=0}^{n-1}\widetilde{j}_{\lambda_j p}e^{-ip\lambda_j}\varepsilon_{n-p}\right)\left(\sum_{q=0}^{n-1}e^{iq\lambda_j}\varepsilon_{n-q}\right).
\end{aligned}
\tag{63}
$$

Because the $\varepsilon_t$ are martingale differences, the second moment of (63) is bounded by

$$
\frac{1}{mn^2}\sum_{j=1}^{m}\sum_{k=1}^{m}\sum_{p=0}^{n-1}|\widetilde{j}_{\lambda_j p}||\widetilde{j}_{\lambda_k p}| + \frac{2}{mn^2}\sum_{j=1}^{m}\sum_{k=1}^{m}\sum_{p=0}^{n-1}|\widetilde{j}_{\lambda_j p}|\sum_{r=0}^{n-1}|\widetilde{j}_{-\lambda_k r}| \tag{64}
$$

$$
+\frac{1}{mn^2}\sum_{j=1}^{m}\sum_{k=1}^{m}\sum_{p=0,p\neq q}^{n-1}|\widetilde{j}_{\lambda_j p}||\widetilde{j}_{-\lambda_k p}|\left|\sum_{q=0}^{n-1}e^{iq(\lambda_j-\lambda_k)}\right|. \tag{65}
$$

Since $\widetilde{j}_{\lambda_j p} = O(\max\{|p|_+^{-1}nj^{-1}, \log n\})$ from Lemma 5.8, (64) is bounded by

$$
\begin{aligned}
&\frac{1}{mn^2}\sum_{j=1}^{m}\sum_{k=1}^{m}\left[\sum_{p=0}^{n-1}(\log n)^2 + \sum_{p=0}^{n-1}\frac{n}{j|p|_+}\sum_{r=0}^{n-1}\frac{n}{k|r|_+}\right] \\
&\quad = O(mn^{-1}(\log n)^2 + m^{-1}(\log n)^4),
\end{aligned}
$$

and, in view of the fact that $\sum_{q=0}^{n-1} e^{iq(\lambda_j-\lambda_k)} = nI\{j=k\}$, (65) is bounded by

$$
\frac{1}{mn}\sum_{j=1}^{m}\sum_{p=0}^{n-1}|\widetilde{j}_{\lambda_j p}|^2 = O\left(\frac{1}{mn}\sum_{j=1}^{m}\sum_{p=0}^{n-1}j^{-1}|p|_+^{-1}n\log n\right) = O(m^{-1}(\log n)^3),
$$

giving (62) $= o_p(1)$. Therefore, we obtain

$$
m^{1/2}\widehat{G}_1(d_0) = 2m^{-1/2}\sum_{j=1}^{m}(\log\lambda_j)I_{uj} + o_p(1).
$$

Let $\nu_j = \log\lambda_j - m^{-1}\sum_1^m \log\lambda_j = \log j - m^{-1}\sum_1^m \log j$ with $\sum_1^m \nu_j = 0$. Then it follows that

$$
\begin{aligned}
m^{1/2}R'(d_0) &= m^{1/2}\left[\frac{\widehat{G}_1(d_0)}{\widehat{G}(d_0)} - 2\frac{1}{m}\sum_{j=1}^{m}\log\lambda_j\right] \\
&= \frac{2m^{-1/2}\sum_1^m(\log\lambda_j)I_{uj} + o_p(1) - (m^{-1}\sum_1^m\log\lambda_j)2m^{-1/2}\sum_1^m I_{uj}}{m^{-1}\sum_1^m I_{uj}} \\
&= \frac{2m^{-1/2}\sum_1^m \nu_j I_{uj} + o_p(1)}{G_0 + o_p(1)}
\end{aligned}
$$

$$= \frac{2m^{-1/2}\sum_1^m \nu_j(I_{uj} - G_0) + o_p(1)}{G_0 + o_p(1)}$$

$$= \frac{2m^{-1/2}\sum_1^m \nu_j(2\pi I_{\varepsilon j} - 1) + o_p(1)}{1 + o_p(1)} \xrightarrow{d} N(0,4),$$

where the fifth line follows from [11], page 1644, completing the proof.

**5. Technical lemmas.** Lemma 5.2 extends Lemma A.3 of Phillips and Shimotsu [8] to hold uniformly in $\theta$. Its proof follows easily from the proof of Lemmas A.2 and A.3 of [8] and is therefore omitted.

LEMMA 5.1 ([7], Theorem 2.2). (a) *If $X_t$ follows* (1), *then*

$$w_u(\lambda) = D_n(e^{i\lambda}; d) w_x(\lambda) - (2\pi n)^{-1/2} e^{in\lambda} \widetilde{X}_{\lambda n}(d),$$

*where* $D_n(e^{i\lambda}; d) = \sum_{k=0}^n \frac{(-d)_k}{k!} e^{ik\lambda}$ *and*

$$\widetilde{X}_{\lambda n}(d) = \widetilde{D}_{n\lambda}(e^{-i\lambda}L; d) X_n = \sum_{p=0}^{n-1} \widetilde{d}_{\lambda p} e^{-ip\lambda} X_{n-p}, \qquad \widetilde{d}_{\lambda p} = \sum_{k=p+1}^{n} \frac{(-d)_k}{k!} e^{ik\lambda}.$$

(b) *If $X_t$ follows* (1) *with $d = 1$, then*

$$w_x(\lambda)(1 - e^{i\lambda}) = w_u(\lambda) - (2\pi n)^{-1/2} e^{i(n+1)\lambda} X_n.$$

LEMMA 5.2 (cf. [8], Lemmas A.2 and A.3). (a) *Uniformly in $\theta \in [-C, C]$ and in $j = 1, 2, \dots, m$ with $m = o(n)$, as $n \to \infty$,*

$$\lambda_j^{-\theta}(1 - e^{i\lambda_j})^\theta = e^{-(\pi/2)\theta i} + O(\lambda_j), \qquad \lambda_j^{-2\theta}|1 - e^{i\lambda_j}|^{2\theta} = 1 + O(\lambda_j^2).$$

(b) *Uniformly in $\theta \in [-1 + \varepsilon, C]$ and in $j = 1, 2, \dots, m$ with $m = o(n)$, as $n \to \infty$,*

$$\lambda_j^{-\theta} D_n(e^{i\lambda_j}; \theta) = e^{-(\pi/2)\theta i} + O(\lambda_j) + O(j^{-1-\theta}),$$
$$\lambda_j^{-2\theta} |D_n(e^{i\lambda_j}; \theta)|^2 = 1 + O(\lambda_j^2) + O(j^{-1-\theta}).$$

LEMMA 5.3. *Let $\widetilde{U}_{\lambda n}(\theta) = \widetilde{D}_{n\lambda}(e^{-i\lambda}L; \theta) u_n = \sum_{p=0}^{n-1} \widetilde{\theta}_{\lambda p} e^{-ip\lambda} u_{n-p}$. Under the assumptions of Theorem* 2.1, *we have, uniformly in $j = 1, \dots, m$, as $n \to \infty$,*

$$E \sup_{\theta \in [-1/2, 1/2]} |n^{\theta - 1/2} j^{1/2 - \theta} \widetilde{U}_{\lambda_j n}(\theta)|^2 = O((\log n)^2).$$

PROOF. When $\theta = 0$, the result follows immediately because $\widetilde{U}_{\lambda_j n}(0) = 0$. When $\theta \neq 0$, define $a_p = \widetilde{\theta}_{\lambda_j p} e^{-ip\lambda_j}$ so that $\widetilde{U}_{\lambda_j n}(\theta) = \sum_{p=0}^{n-1} a_p u_{n-p}$. We suppress the dependence of $a_p$ on $\theta$ and $\lambda_j$. Summation by parts gives

$$\widetilde{U}_{\lambda_j n}(\theta) = \sum_{p=0}^{n-2}(a_p - a_{p+1}) \sum_{q=0}^{p} u_{n-q} + a_{n-1} \sum_{q=0}^{n-1} u_{n-q}.$$

Phillips and Shimotsu ([8], page 670) show that (note that Phillips and Shimotsu use $\lambda_s$ instead of $\lambda_j$ to denote Fourier frequencies)

$$a_p - a_{p+1} = b_{np}(\theta) + \frac{(-\theta)_n}{n!} e^{-ip\lambda_j},$$

where

$$b_{np}(\theta) = \sum_{k=p+1}^{n-1} \frac{(1+\theta)\Gamma(k-\theta)}{\Gamma(-\theta)\Gamma(k+2)} e^{i(k-p)\lambda_j}. \tag{66}$$

Then, since $a_{n-1} = (-\theta)_n e^{-i(n-1)\lambda_j}/n!$, we have

$$\begin{aligned}
\widetilde{U}_{\lambda_j n}(\theta) &= \sum_{p=0}^{n-2} b_{np}(\theta) \sum_{q=0}^{p} u_{n-q} \\
&\quad + \frac{(-\theta)_n}{n!} \sum_{p=0}^{n-2} e^{-ip\lambda_j} \sum_{q=0}^{p} u_{n-q} + \frac{(-\theta)_n}{n!} e^{-i(n-1)\lambda_j} \sum_{q=0}^{n-1} u_{n-q} \\
&= \sum_{p=0}^{n-2} b_{np}(\theta) \sum_{q=0}^{p} u_{n-q} + \frac{(-\theta)_n}{n!} \sum_{p=0}^{n-1} e^{-ip\lambda_j} \sum_{q=0}^{p} u_{n-q} \\
&= U_{1n}(\theta) + U_{2n}(\theta).
\end{aligned} \tag{67}$$

We proceed to show that the elements of $n^{\theta-1/2} j^{1/2-\theta} U_{\cdot n}(\theta)$ are of the stated order. First, for $U_{1n}$, we have

$$\sup_{\theta} |n^{\theta-1/2} j^{1/2-\theta} U_{1n}(\theta)| \leq \sum_{p=0}^{n-2} \sup_{\theta} |n^{\theta-1/2} j^{1/2-\theta} b_{np}(\theta)| \left| \sum_{q=0}^{p} u_{n-q} \right|.$$

Because $\sum_{-\infty}^{\infty} E u_t u_{t+q} = 2\pi f_u(0) = 2\pi G_0 < \infty$, it follows from Kronecker's lemma that, uniformly in $p = 0, \ldots, n-1$,

$$E\left(\sum_{q=0}^{p} u_{n-q}\right)^2 = (p+1) \sum_{q=-p}^{p} (1 - |q|/(p+1)) E u_t u_{t+q} = O(|p|_+). \tag{68}$$

Therefore, if we have, uniformly in $p = 0, \ldots, n-1$ and $j = 1, \ldots, m$,

$$\sup_{\theta \in [-1/2, 1/2]} |n^{\theta-1/2} j^{1/2-\theta} b_{np}(\theta)| = O(|p|_+^{-3/2}), \tag{69}$$

it follows from Minkowski's inequality that

$$E\sup_{\theta}|n^{\theta-1/2}j^{1/2-\theta}U_{1n}(\theta)|^2 = O\left(\left(\sum_{p=0}^{n-2}|p|_+^{-1}\right)^2\right) = O((\log n)^2). \tag{70}$$

To show (69), Phillips and Shimotsu ([8], page 670, equation (21)) show that

$$|b_{np}(\theta)| = O(\min\{|p|_+^{-\theta-1}, |p|_+^{-\theta-2}nj^{-1}\}) \tag{71}$$

uniformly in $\theta \in [-\frac{1}{2}, \frac{1}{2}]$, $p = 0, \ldots, n-1$, and $j = 1, \ldots, m$. Although Phillips and Shimotsu do not state explicitly that the bound (71) holds uniformly in $\theta \in [-\frac{1}{2}, \frac{1}{2}]$, it is clear from its proof that (71) holds uniformly in $\theta \in [-\frac{1}{2}, \frac{1}{2}]$. Then (69) follows from (71) because

$$0 \le p \le n/j : n^{\theta-1/2}j^{1/2-\theta}|p|_+^{-\theta-1} = (j|p|_+/n)^{1/2-\theta}|p|_+^{-3/2} \le |p|_+^{-3/2},$$

$$n/j \le p \le n : n^{\theta-1/2}j^{1/2-\theta}p^{-\theta-2}nj^{-1} = (jp/n)^{-\theta-1/2}p^{-3/2} \le p^{-3/2}.$$

For $U_{2n} = ((-\theta)_n/n!)\sum_0^{n-1} e^{-ip\lambda_j}\sum_0^p u_{n-q}$, first we rewrite the sum as

$$\begin{aligned}
\sum_{p=0}^{n-1} e^{-ip\lambda_j}\sum_{q=0}^{p} u_{n-q} &= \sum_{n-p=1}^{n} e^{i(n-p)\lambda_j}\sum_{n-q=n-p}^{n} u_{n-q} \\
&= \sum_{k=1}^{n} u_k \sum_{q=1}^{k} e^{iq\lambda_j} \\
&= \sum_{k=1}^{n} u_k \frac{e^{i\lambda_j}(1-e^{ik\lambda_j})}{1-e^{i\lambda_j}} \\
&= \frac{e^{i\lambda_j}}{1-e^{i\lambda_j}}\sum_{k=1}^{n} u_k - \frac{e^{i\lambda_j}}{1-e^{i\lambda_j}}(2\pi n)^{1/2} w_u(\lambda_j).
\end{aligned} \tag{72}$$

Since $(-\theta)_n/n! = O(n^{-\theta-1})$ uniformly in $\theta \in [-\frac{1}{2}, \frac{1}{2}]$ and $(1-e^{i\lambda_j})^{-1} = O(nj^{-1})$, $E\sup_\theta|n^{\theta-1/2}j^{1/2-\theta}U_{2n}|^2 = O(1)$ follows from (68) and $E|w_u(\lambda_j)|^2 = O(1)$ ([11], page 1637). □

LEMMA 5.4. *For $\kappa \in (0,1)$ and $C \in (1,\infty)$, as $m \to \infty$,*

(a) $$\sup_{-C\le\gamma\le C}\left|\frac{1}{m}\sum_{j=[\kappa m]}^{m}\left(\frac{j}{m}\right)^\gamma - \int_\kappa^1 x^\gamma\,dx\right| = O(m^{-1}),$$

(b) $$\sup_{-C\le\gamma\le C}\left|m^{-1}\sum_{j=[\kappa m]}^{m}(j/m)^\gamma\right| = O(1),$$

$$\liminf_{m\to\infty}\inf_{-C\le\gamma\le C}\left(m^{-1}\sum_{j=[\kappa m]}^{m}(j/m)^\gamma\right) > \varepsilon > 0.$$

PROOF. Note that $[\kappa m] \geq 3$ for large $m$. For part (a), since

$$\frac{1}{m}\sum_{j=[\kappa m]}^{m}\left(\frac{j}{m}\right)^{\gamma} = \sum_{j=[\kappa m]}^{m}\int_{(j-1)/m}^{j/m}\left(\frac{j}{m}\right)^{\gamma} dx,$$

$$\int_{\kappa}^{1} x^{\gamma}\, dx = \sum_{j=[\kappa m]}^{m}\int_{(j-1)/m}^{j/m} x^{\gamma}\, dx - \int_{([\kappa m]-1)/m}^{\kappa} x^{\gamma}\, dx,$$

their difference is bounded uniformly in $\gamma$ by, for sufficiently large $m$:

$$\sum_{j=[\kappa m]}^{m}\left|\int_{(j-1)/m}^{j/m}\left\{\left(\frac{j}{m}\right)^{\gamma} - x^{\gamma}\right\} dx\right| + \int_{\kappa-(2/m)}^{\kappa} x^{\gamma}\, dx$$

$$\leq \frac{c}{m^2}\sum_{j=[\kappa m]}^{m}\left(\frac{j}{m}\right)^{-C-1} + \frac{c}{m} = O(m^{-1}),$$

by the mean value theorem. Part (b) follows immediately from part (a). □

LEMMA 5.5. *For $p \sim m/e$ as $m \to \infty$ and $\Delta \in (0, \frac{1}{2e})$, there exist $\varepsilon \in (0, 0.1)$ and $\bar{\kappa} \in (0, \frac{1}{4})$ such that, for all fixed $\kappa \in (0, \bar{\kappa}]$ and sufficiently large $m$:*

(a) $$\inf_{-C \leq \gamma \leq -1+2\Delta}\frac{1}{m}\sum_{j=[\kappa m]}^{m}\left(\frac{j}{p}\right)^{\gamma} \geq 1 + 2\varepsilon,$$

(b) $$\inf_{1 \leq \gamma \leq C}\frac{1}{m}\sum_{j=[\kappa m]}^{m}\left(\frac{j}{p}\right)^{\gamma} \geq 1 + 2\varepsilon.$$

PROOF. From Lemma 5.4 we obtain, for large $m$,

$$\begin{aligned}\inf_{-C \leq \gamma \leq -1+2\Delta}\frac{1}{m}\sum_{j=[\kappa m]}^{m}\left(\frac{j}{p}\right)^{\gamma} &\geq \inf_{-C \leq \gamma \leq -1+2\Delta}\frac{1}{m}\sum_{j=[\kappa m]}^{p}\left(\frac{j}{p}\right)^{\gamma} \\ &\geq \frac{1}{m}\sum_{j=[\kappa m]}^{p}\left(\frac{j}{p}\right)^{-1+2\Delta} \\ &\sim \frac{1}{e}\int_{\kappa e}^{1} x^{2\Delta-1}\, dx = \frac{1-(\kappa e)^{2\Delta}}{2\Delta e},\end{aligned}$$

$$\inf_{1 \leq \gamma \leq C}\frac{1}{m}\sum_{j=[\kappa m]}^{m}\left(\frac{j}{p}\right)^{\gamma} \sim \inf_{1 \leq \gamma \leq C}\frac{e^{\gamma}}{\gamma+1}(1-\kappa^{\gamma+1}) \geq \frac{e}{2}(1-\kappa^2),$$

where the last inequality holds because $e^{\gamma}/(\gamma+1)$ is monotone increasing for $\gamma \geq 1$. Since $2\Delta e < 1$, choosing $\kappa$ sufficiently small gives the stated results. □

Lemma 5.6. *For $\kappa \in (0,1)$, $C \in (1,\infty)$ and $m = o(n)$, as $n \to \infty$,*

$$E \sup_{\alpha \in [-C,C]} \left| \frac{1}{m} \sum_{j=[\kappa m]}^{m} \left(\frac{j}{m}\right)^{\alpha} w_u(\lambda_j) \right| = O(m^{-1/2} \log m).$$

Proof. Summation by parts gives

$$\begin{aligned} &\frac{1}{m} \sum_{j=[\kappa m]}^{m} \left(\frac{j}{m}\right)^{\alpha} w_u(\lambda_j) \\ &\quad = \frac{1}{m} \sum_{r=[\kappa m]}^{m-1} \left[\left(\frac{r}{m}\right)^{\alpha} - \left(\frac{r+1}{m}\right)^{\alpha}\right] \sum_{j=[\kappa m]}^{r} w_u(\lambda_j) + \frac{1}{m} \sum_{j=[\kappa m]}^{m} w_u(\lambda_j). \end{aligned}$$

Note that, uniformly in $r = 1, \dots, m-1$ and $\alpha$,

$$\left|\left(\frac{r}{m}\right)^{\alpha} - \left(\frac{r+1}{m}\right)^{\alpha}\right| = \left(\frac{r}{m}\right)^{\alpha} \left|1 - \left(1 + \frac{1}{r}\right)^{\alpha}\right| \le c \left(\frac{r}{m}\right)^{-C} \frac{1}{r}, \tag{73}$$

because $\sup_{\alpha} |(1+x)^{\alpha} - 1| \le C 2^C x$ for $0 \le x \le 1$. The results in ([11], page 1637) imply that $E|w_u(\lambda_j) - C(e^{i\lambda_j}) w_\varepsilon(\lambda_j)|^2 = O(j^{-1} \log(j+1))$ uniformly in $j = 1, \dots, m$, giving

$$E \left| \sum_{j=[\kappa m]}^{r} w_u(\lambda_j) \right|^2 = O(r \log(r+1)), \qquad r = [\kappa m], \dots, m. \tag{74}$$

From (73), (74) and Lemma 5.4, $E \sup_{\alpha \in [-C,C]} |m^{-1} \sum_{[\kappa m]}^{m} (j/m)^{\alpha} w_u(\lambda_j)|$ is bounded by

$$m^{-1} \sum_{r=[\kappa m]}^{m-1} (r/m)^{-C} r^{-1/2} \log r + m^{-1/2} \log m = O(m^{-1/2} \log m),$$

giving the required result. □

Lemma 5.7. *Define $J_n(L) = \sum_{k=1}^{n} L^k / k$ and $D_n(L; d) = \sum_{k=0}^{n} \frac{(-d)_k}{k!} L^k$. Then:*

(a) $$J_n(L) = J_n(e^{i\lambda}) + \widetilde{J}_{n\lambda}(e^{-i\lambda} L)(e^{-i\lambda} L - 1),$$

(b) $$\begin{aligned} J_n(L) D_n(L; d) &= J_n(e^{i\lambda}) D_n(e^{i\lambda}; d) + D_n(e^{i\lambda}; d) \widetilde{J}_{n\lambda}(e^{-i\lambda} L)(e^{-i\lambda} L - 1) \\ &\quad + J_n(L) \widetilde{D}_{n\lambda}(e^{-i\lambda} L; d)(e^{-i\lambda} L - 1), \end{aligned}$$

*where*

$$\widetilde{J}_{n\lambda}(e^{-i\lambda} L) = \sum_{p=0}^{n-1} \widetilde{j}_{\lambda p} e^{-ip\lambda} L^p, \qquad \widetilde{j}_{\lambda p} = \sum_{p+1}^{n} \frac{1}{k} e^{ik\lambda},$$

$$\widetilde{D}_{n\lambda}(e^{-i\lambda} L; d) = \sum_{p=0}^{n-1} \widetilde{d}_{\lambda p} e^{-ip\lambda} L^p, \qquad \widetilde{d}_{\lambda p} = \sum_{p+1}^{n} \frac{(-d)_k}{k!} e^{ik\lambda}.$$

PROOF. For part (a) see [9], formula (32). For part (b), from Lemma 2.1 of [7] we have $D_n(L;d) = D_n(e^{i\lambda};d) + \widetilde{D}_{n\lambda}(e^{-i\lambda}L;d)(e^{-i\lambda}L - 1)$, and the stated result follows immediately. □

LEMMA 5.8. *Let* $J_n(e^{i\lambda}) = \sum_{k=1}^n e^{ik\lambda}/k$ *and* $\widetilde{j}_{\lambda p} = \sum_{k=p+1}^n e^{ik\lambda}/k$, *as defined in Lemma 5.7. Then uniformly in* $p = 0, \ldots, n-1$ *and* $j = 1, \ldots, m$ *with* $m = o(n)$, *as* $n \to \infty$:

(a) $J_n(e^{i\lambda_j}) = -\log\lambda_j + \frac{i}{2}(\pi - \lambda_j) + O(\lambda_j^2) + O(j^{-1})$,

(b) $\widetilde{j}_{\lambda_j p} = O(\min\{|p|_+^{-1} n j^{-1}, \log n\})$.

PROOF. For (a), first we have

$$J_n(e^{i\lambda_j}) = \sum_{k=1}^n \frac{1}{k} e^{ik\lambda_j} = \sum_{k=1}^\infty \frac{1}{k} e^{ik\lambda_j} - \sum_{k=n+1}^\infty \frac{1}{k} e^{ik\lambda_j}. \tag{75}$$

The first term on the right-hand side of (75) is equal to ([16], page 5)

$$\sum_{k=1}^\infty \frac{\cos k\lambda_j}{k} + i\sum_{k=1}^\infty \frac{\sin k\lambda_j}{k} = -\log\left|2\sin\frac{\lambda_j}{2}\right| + i\frac{1}{2}(\pi - \lambda_j).$$

Since $2\sin(\lambda_j/2) = \lambda_j + O(\lambda_j^3) = \lambda_j(1 + O(\lambda_j^2))$, the right-hand side is equal to

$$\begin{aligned} &-\log\lambda_j - \log(1 + O(\lambda_j^2)) + i\frac{1}{2}(\pi - \lambda_j) \\ &\quad = -\log\lambda_j + O(\lambda_j^2) + \frac{i}{2}(\pi - \lambda_j). \end{aligned}$$

For the second term on the right-hand side of (75), summation by parts gives

$$\begin{aligned} &\left|\sum_{k=n+1}^\infty \frac{1}{k} e^{ik\lambda_j}\right| \\ &\quad = \left|\lim_{N\to\infty}\left[\sum_{r=n+1}^{n+N-1}\left(\frac{1}{r} - \frac{1}{r+1}\right)\sum_{k=n+1}^r e^{ik\lambda_j} + \frac{1}{n+N}\sum_{k=n+1}^{n+N} e^{ik\lambda_j}\right]\right| \\ &\quad \le C\left(nj^{-1}\sum_{k=n}^\infty r^{-2} + j^{-1}\right) = O(j^{-1}), \end{aligned}$$

giving (a). Part (b) follows from the fact that

$$|\widetilde{j}_{\lambda_j p}| \le (p+1)^{-1} \max_{p+1 \le N \le n}\left|\sum_{k=p+1}^N e^{ik\lambda_j}\right| \quad \text{and} \quad \widetilde{j}_{\lambda_j p} = O\left(\sum_{k=0}^n |k|_+^{-1}\right). \qquad \square$$

LEMMA 5.9. *Suppose* $Y_t = (1-L)^\theta u_t I\{t \geq 1\}$ . *Under the assumptions of Theorem 2.2 we have:*

(a) $-w_{\log(1-L)y}(\lambda_j) = J_n(e^{i\lambda_j}) D_n(e^{i\lambda_j};\theta) w_u(\lambda_j) + n^{-1/2} V_{nj}(\theta)$,

(b) $-w_{\log(1-L)u}(\lambda_j) = J_n(e^{i\lambda_j}) w_u(\lambda_j)$
$$- C(1)(2\pi n)^{-1/2} \widetilde{J}_{n\lambda_j}(e^{-i\lambda_j}L)\varepsilon_n + r_{nj},$$

(c) $w_{(\log(1-L))^2 y}(\lambda_j) = J_n(e^{i\lambda_j})^2 D_n(e^{i\lambda_j};\theta) w_u(\lambda_j) + n^{-1/2}\Psi_{nj}(\theta)$,

*where, uniformly in* $j = 1, \dots, m$, *as* $n \to \infty$,

$$E \sup_\theta |n^{\theta-1/2} j^{1/2-\theta} V_{nj}(\theta)|^2 = O((\log n)^4),$$

$$E|j^{1/2} r_{nj}|^2 = o(1) + O(j^{-1}),$$

$$E \sup_\theta |n^{\theta-1/2} j^{1/2-\theta} \Psi_{nj}(\theta)|^2 = O((\log n)^6).$$

PROOF. Define $\bar{u}_t = u_t I\{t \geq 1\}$, so that $Y_t = D_{t-1}(L;\theta)\bar{u}_t = D_n(L;\theta)\bar{u}_t$ for $t \leq n$. Since $Y_t = 0$ for $t \leq 0$, we have

$$\log(1-L)Y_t = (-L - L^2/2 - L^3/3 - \cdots)Y_t = -J_n(L)Y_t.$$

For parts (a) and (b), from Lemma 5.7(b) we have

$$\begin{aligned} &-\log(1-L)Y_t \\ &\quad = J_n(L)D_n(L;\theta)\bar{u}_t \\ &\quad = J_n(e^{i\lambda_j})D_n(e^{i\lambda_j};\theta)\bar{u}_t + D_n(e^{i\lambda_j};\theta)\widetilde{J}_{n\lambda_j}(e^{-i\lambda_j}L)(e^{-i\lambda_j}L-1)\bar{u}_t \\ &\quad\quad + J_n(L)\widetilde{D}_{n\lambda_j}(e^{-i\lambda_j}L;\theta)(e^{-i\lambda_j}L-1)\bar{u}_t. \end{aligned}$$

Since $\sum_{t=1}^n e^{it\lambda_j}(e^{-i\lambda_j}L-1)\bar{u}_t = -\bar{u}_n$, taking the d.f.t. of the right-hand side gives

$$\begin{aligned} &J_n(e^{i\lambda_j})D_n(e^{i\lambda_j};\theta)w_u(\lambda_j) - (2\pi n)^{-1/2} D_n(e^{i\lambda_j};\theta)\widetilde{J}_{n\lambda_j}(e^{-i\lambda_j}L)\bar{u}_n \\ &\quad - (2\pi n)^{-1/2} J_n(L)\widetilde{D}_{n\lambda_j}(e^{-i\lambda_j}L;\theta)\bar{u}_n. \end{aligned} \tag{76}$$

Note that Lemma 5.2(b) gives $|D_n(e^{i\lambda_j};\theta)| \leq c\lambda_j^\theta$. Therefore part (a) follows if

$$E|\widetilde{J}_{n\lambda_j}(e^{-i\lambda_j}L)\bar{u}_n|^2 = O(nj^{-1}), \tag{77}$$

$$E \sup_\theta |n^{\theta-1/2} j^{1/2-\theta} J_n(L)\widetilde{D}_{n\lambda_j}(e^{-i\lambda_j}L;\theta)\bar{u}_n|^2 = O((\log n)^4). \tag{78}$$

First we show (77). Define $a_p' = \widetilde{j}_{\lambda_j p} e^{-ip\lambda_j} = \sum_{k=p+1}^n k^{-1} e^{i(k-p)\lambda_j}$, so that $\widetilde{J}_{n\lambda_j}(e^{-i\lambda_j}L)\bar{u}_n = \sum_{p=0}^{n-1} a_p' \bar{u}_{n-p} = \sum_{p=0}^{n-1} a_p' u_{n-p}$. Then summation by parts

gives

$$\widetilde{J}_{n\lambda_j}(e^{-i\lambda_j}L)\bar{u}_n = \sum_{p=0}^{n-2}(a'_p - a'_{p+1})\sum_{q=0}^{p} u_{n-q} + a'_{n-1}\sum_{q=0}^{n-1} u_{n-q}.$$

Observe that

$$\begin{aligned} a'_p - a'_{p+1} &= \sum_{k=p+1}^{n} \frac{1}{k} e^{i(k-p)\lambda_j} - \sum_{k=p+2}^{n} \frac{1}{k} e^{i(k-p-1)\lambda_j} \\ &= \sum_{k=p+1}^{n-1}\left[\frac{1}{k} - \frac{1}{k+1}\right] e^{i(k-p)\lambda_j} + \frac{1}{n} e^{-ip\lambda_j} \\ &= \sum_{k=p+1}^{n-1} \frac{1}{k(k+1)} e^{i(k-p)\lambda_j} + \frac{1}{n} e^{-ip\lambda_j}. \end{aligned}$$

Define $c_{np} = \sum_{k=p+1}^{n-1} \frac{1}{k(k+1)} e^{i(k-p)\lambda_j}$. Then since $a'_{n-1} = n^{-1}e^{-i(n-1)\lambda_j}$, we obtain

$$\begin{aligned} &\widetilde{J}_{n\lambda_j}(e^{-i\lambda_j}L)\bar{u}_n \\ &\quad = \sum_{p=0}^{n-2} c_{np}\sum_{q=0}^{p} u_{n-q} + \frac{1}{n}\sum_{p=0}^{n-2} e^{-ip\lambda_j}\sum_{q=0}^{p} u_{n-q} + \frac{1}{n} e^{-i(n-1)\lambda_j}\sum_{q=0}^{n-1} u_{n-q} \\ &\quad = \sum_{p=0}^{n-2} c_{np}\sum_{q=0}^{p} u_{n-q} + \frac{1}{n}\sum_{p=0}^{n-1} e^{-ip\lambda_j}\sum_{q=0}^{p} u_{n-q} \qquad (79) \\ &\quad = \sum_{p=0}^{n-2} c_{np}\sum_{q=0}^{p} u_{n-q} + \left[\frac{1}{n}\frac{e^{i\lambda_j}}{1-e^{i\lambda_j}}\sum_{k=1}^{n} u_k - \frac{1}{n}\frac{e^{i\lambda_j}}{1-e^{i\lambda_j}}(2\pi n)^{1/2} w_u(\lambda_j)\right] \\ &\quad = \widetilde{J}_{1n} + \widetilde{J}_{2n}, \end{aligned}$$

where the fourth line follows from (72). $E|\widetilde{J}_{2n}|^2 = O(nj^{-2})$ in view of the order of magnitude of $E|\sum_1^n u_k|^2$ and $E|w_u(\lambda_j)|^2$. For $\widetilde{J}_{1n}$, since

$$\begin{aligned} |c_{np}| &= \left|\sum_{k=p+1}^{n-1} \frac{1}{k(k+1)} e^{i(k-p)\lambda_j}\right| \\ &\le |p|_+^{-2} \max_{1\le N\le n}\left|\sum_{p+1}^{p+N} e^{ik\lambda_j}\right| \le C|p|_+^{-2} nj^{-1}, \end{aligned}$$

$$|c_{np}| \le \left|\sum_{k=p+1}^{n-1} \frac{1}{k(k+1)}\right| \le C|p|_+^{-1},$$

we have

$$|c_{np}| \le C \min\{|p|_+^{-1}, |p|_+^{-2} n j^{-1}\}. \tag{80}$$

Therefore, it follows from (68) and Minkowski's inequality that

$$E|\widetilde{J}_{1n}|^2 = O\left(\left(\sum_{p=0}^{n/j} |p|_+^{-1/2} + \sum_{p=n/j}^{n} \frac{n}{j} |p|_+^{-3/2}\right)^2\right) = O(nj^{-1}), \tag{81}$$

and hence (77) follows.

Now we move to the proof of (78). When $\theta = 0$, then $\widetilde{D}_{n\lambda_j}(e^{-i\lambda_j} L; \theta) = 0$, and (78) follows immediately. Assume $\theta \neq 0$. If we have, uniformly in $r = 0, 1, \dots,$

$$E \sup_\theta |n^{\theta-1/2} j^{1/2-\theta} L^r \widetilde{D}_{n\lambda_j}(e^{-i\lambda_j} L; \theta) \bar{u}_n|^2 = O((\log n)^2), \tag{82}$$

then (78) follows because Minkowski's inequality gives

$$\begin{aligned} E \sup_\theta &|n^{\theta-1/2} j^{1/2-\theta} J_n(L) \widetilde{D}_{n\lambda_j}(e^{-i\lambda_j} L; \theta) \bar{u}_n|^2 \\ &\le E\left(\sum_{p=1}^{n-1} p^{-1} \sup_\theta |n^{\theta-1/2} j^{1/2-\theta} L^p \widetilde{D}_{n\lambda_j}(e^{-i\lambda_j} L; \theta) \bar{u}_n|\right)^2 \\ &\le \left(\sum_{p=1}^{n-1} p^{-1} \left(E \sup_\theta |n^{\theta-1/2} j^{1/2-\theta} L^p \widetilde{D}_{n\lambda_j}(e^{-i\lambda_j} L; \theta) \bar{u}_n|^2\right)^{1/2}\right)^2 \\ &= O((\log n)^4). \end{aligned}$$

We proceed to show (82). For $r \ge n$, (82) follows immediately because $L^r \widetilde{D}_{n\lambda_j}(e^{-i\lambda_j} L; \theta) \bar{u}_n = 0$. For $r = 0, \dots, n-1$, using a decomposition similar to (67) gives

$$\begin{aligned} L^r &\widetilde{D}_{n\lambda_j}(e^{-i\lambda_j} L; \theta) \bar{u}_n \\ &= \sum_{p=0}^{n-2} b_{np}(\theta) L^r \sum_{q=0}^{p} \bar{u}_{n-q} + \frac{(-\theta)_n}{n!} L^r \sum_{p=0}^{n-1} e^{-ip\lambda_j} \sum_{q=0}^{p} \bar{u}_{n-q} \\ &= U'_{1n}(\theta) + U'_{2n}(\theta), \end{aligned}$$

where $b_{np}(\theta)$ is defined in (66). For $U'_{1n}(\theta)$, since $E(L^r \sum_{q=0}^{p} \bar{u}_{n-q})^2 = O(|p|_+^{1/2})$, the arguments in the proof of Lemma 5.3 go through and $E \sup_\theta |n^{\theta-1/2} j^{1/2-\theta} \times U'_{1n}(\theta)|^2 = O((\log n)^2)$ holds. For $U'_{2n}(\theta)$, using a decomposition similar to (72) gives

$$U'_{2n}(\theta) = \frac{(-\theta)_n}{n!} \frac{e^{i\lambda_j}}{1 - e^{i\lambda_j}} L^r \sum_{k=1}^{n} \bar{u}_k - \frac{(-\theta)_n}{n!} \frac{e^{i\lambda_j}}{1 - e^{i\lambda_j}} L^r \sum_{k=1}^{n} e^{ik\lambda_j} \bar{u}_k$$

$$= \frac{(-\theta)_n}{n!} \frac{e^{i\lambda_j}}{1-e^{i\lambda_j}} \sum_{k=1}^{n-r} u_k - \frac{(-\theta)_n}{n!} \frac{e^{i\lambda_j}}{1-e^{i\lambda_j}} e^{ir\lambda_j} \sum_{q=1}^{n-r} e^{iq\lambda_j} u_q.$$

Since $E(\sum_{k=1}^{n-r} u_k)^2 = O(n^{1/2})$ for any $r$, $E\sup_\theta |n^{\theta-1/2} j^{1/2-\theta} U'_{2n}(\theta)|^2 = O(1)$ and (82) follow if, for $m = o(n)$,

$$\max_{1\le r\le n} \max_{1\le j\le m} E\left|(2\pi r)^{-1/2} \sum_{k=1}^{r} e^{ik\lambda_j} u_k\right|^2 = O(1). \tag{83}$$

We establish (83) to complete the proof of part (a). An elementary calculation gives

$$E\left|(2\pi r)^{-1/2} \sum_{k=1}^{r} e^{ik\lambda_j} u_k\right|^2 = \int_{-\pi}^{\pi} f_u(\lambda) K_r(\lambda-\lambda_j)\, d\lambda,$$

where $K_r(\lambda) = (2\pi r)^{-1} \sum_{s=1}^{r} \sum_{t=1}^{r} e^{i(t-s)\lambda}$ is Fejér's kernel. From Zygmund [16], pages 88–90, $\int_{-\pi}^{\pi} |K_r(\lambda)|\, d\lambda < A$ and $|K_r(\lambda)| \le Ar^{-1}\lambda^{-2}$ for a finite constant $A$. Furthermore, from Assumption 1 there exists $\eta \in (0,\pi)$ such that $\sup_{\lambda\in[-\eta,\eta]} |f_u(\lambda)| < C$, and $\inf_{|\lambda|>\eta} |\lambda-\lambda_j| \ge \eta/2$ if $\lambda_j < \eta/2$. It follows that for sufficiently large $n$

$$\begin{aligned}
&\int_{-\pi}^{\pi} f_u(\lambda) K_r(\lambda-\lambda_j)\, d\lambda \\
&\quad = \int_{|\lambda|\le\eta} f_u(\lambda) K_r(\lambda-\lambda_j)\, d\lambda + \int_{\eta\le|\lambda|\le\pi} f_u(\lambda) K_r(\lambda-\lambda_j)\, d\lambda \\
&\quad \le AC + Ar^{-1}(\eta/2)^{-2} \int_{\eta\le|\lambda|\le\pi} f_u(\lambda)\, d\lambda < \infty,
\end{aligned}$$

uniformly in $j = 1,\dots,m$, and (83) follows.

For part (b), in view of (76), $D_n(e^{i\lambda_j};0) = 1$ and $\widetilde{D}_{n\lambda_j}(e^{-i\lambda_j}L;0) = 0$, part (b) follows if, as $n\to\infty$, uniformly in $j = 1,\dots,m$,

$$E|j^{1/2} n^{-1/2} \widetilde{J}_{n\lambda_j}(e^{-i\lambda_j}L)(\bar{u}_n - C(1)\varepsilon_n)|^2 = o(1) + O(j^{-1}). \tag{84}$$

Using the same decomposition as (79), write $j^{1/2} n^{-1/2} \widetilde{J}_{n\lambda_j}(e^{-i\lambda_j}L)(\bar{u}_n - C(1)\varepsilon_n)$ as

$$\sum_{p=0}^{n-2} \frac{j^{1/2}}{\sqrt{n}} c_{np} \sum_{q=0}^{p} (u_{n-q} - C(1)\varepsilon_{n-q}) \tag{85}$$

$$\begin{aligned}
&+ \frac{j^{1/2}}{n\sqrt{n}} \frac{e^{i\lambda_j}}{1-e^{i\lambda_j}} \sum_{k=1}^{n} (u_{n-k} - C(1)\varepsilon_{n-k}) \\
&- \frac{j^{1/2}}{n} \frac{e^{i\lambda_j}\sqrt{2\pi}}{1-e^{i\lambda_j}} [w_u(\lambda_j) - C(1) w_\varepsilon(\lambda_j)].
\end{aligned} \tag{86}$$

If we have

$$E\left[\sum_{q=0}^{p}(u_{n-q}-C(1)\varepsilon_{n-q})\right]^2 = \begin{cases} O(|p|_+), & \text{uniformly in } p=0,\dots,n-1, \\ o(p), & \text{as } p\to\infty, \end{cases} \tag{87}$$

then it follows from Minkowski's inequality and the order of $c_{np}$ given by (80) that

$$\begin{aligned}(E|(85)|^2)^{1/2} &= O\left((j/n)^{1/2}\sum_{p=0}^{\sqrt{n/j}}|p|_+^{-1/2}\right) \\ &\quad + o\left((j/n)^{1/2}\sum_{p=\sqrt{n/j}}^{n/j}p^{-1/2}+(j/n)^{1/2}\sum_{p=n/j}^{n}\frac{n}{j}p^{-3/2}\right) \\ &= O((j/n)^{1/4})+o(1) \\ &= o(1),\end{aligned}$$

because $\sqrt{n/j}\ge\sqrt{n/m}\to\infty$ from Assumption $4'$. To prove (87), note that when $p=0$, (87) follows immediately. When $p\ge 1$, observe that

$$\begin{aligned}&E\left[\sum_{q=0}^{p}(u_{n-q}-C(1)\varepsilon_{n-q})\right]^2 \\ &\qquad \le 2E\left[\sum_{q=0}^{p}u_{n-q}\right]^2+2E\left[C(1)\sum_{q=0}^{p}\varepsilon_{n-q}\right]^2.\end{aligned}$$

Since the first term on the right-hand side is uniformly $O(p)$ from (68) and the second term on the right-hand side is equal to $2C(1)^2(p+1)$, the first part of (87) holds. For the second part of (87), note that the left-hand side of (87) is equal to ($\gamma_q=Eu_tu_{t+q}$)

$$\begin{aligned}&\sum_{r=-p}^{p}(p+1-|r|)\gamma_r-2C(1)\sum_{q=0}^{p}\sum_{r=0}^{q}c_{q-r}+(p+1)C(1)^2 \\ &\quad=-(p+1)\sum_{|r|\ge p+1}\gamma_r-2\sum_{r=1}^{p}r\gamma_r+2C(1)(p+1)\sum_{r\ge p+1}c_r-2C(1)\sum_{r=1}^{p}rc_r.\end{aligned}$$

If $\sum_{-\infty}^{\infty}a_r$ converges, then $\sum_{|r|\ge p+1}a_r$ tends to 0 as $p\to\infty$; thus the first and third terms are $o(p)$ because both $\sum_{-\infty}^{\infty}\gamma_r$ and $\sum_{-\infty}^{\infty}c_r$ converge. The second and fourth terms are $o(p)$ from Kronecker's lemma, and the second part of (87) follows. Obviously $E|(86)|^2=O(j^{-1})$, and (84) follows.

For part (c), first from Lemma 2.1 of [7] and Lemma 5.7 we have

$$\begin{aligned}J_n(L)^2 &= J_n(L)[J_n(e^{i\lambda}) + \widetilde{J}_{n\lambda}(e^{-i\lambda}L)(e^{-i\lambda}L-1)]\\ &= J_n(L)J_n(e^{i\lambda}) + J_n(L)\widetilde{J}_{n\lambda}(e^{-i\lambda}L)(e^{-i\lambda}L-1)\\ &= J_n(e^{i\lambda})^2 + J_n(e^{i\lambda})\widetilde{J}_{n\lambda}(e^{-i\lambda}L)(e^{-i\lambda}L-1)\\ &\quad + J_n(L)\widetilde{J}_{n\lambda}(e^{-i\lambda}L)(e^{-i\lambda}L-1),\\ D_n(L;\theta) &= D_n(e^{i\lambda};\theta) + \widetilde{D}_{n\lambda}(e^{-i\lambda}L;\theta)(e^{-i\lambda}L-1).\end{aligned}$$

It follows that

$$\begin{aligned}(\log(1-L))^2 Y_t &= J_n(L)^2 D_n(L;\theta)\bar{u}_t\\ &= J_n(e^{i\lambda})^2 D_n(e^{i\lambda};\theta)\bar{u}_t\\ &\quad + D_n(e^{i\lambda};\theta)[J_n(e^{i\lambda}) + J_n(L)]\widetilde{J}_{n\lambda}(e^{-i\lambda}L)(e^{-i\lambda}L-1)\bar{u}_t\\ &\quad + J_n(L)^2\widetilde{D}_{n\lambda}(e^{-i\lambda}L;\theta)(e^{-i\lambda}L-1)\bar{u}_t.\end{aligned}$$

Taking its d.f.t. gives

$$\begin{aligned}&J_n(e^{i\lambda_j})^2 D_n(e^{i\lambda_j};\theta)w_u(\lambda_j)\\ &\quad - (2\pi n)^{-1/2} D_n(e^{i\lambda_j};\theta)[J_n(e^{i\lambda_j}) + J_n(L)]\widetilde{J}_{n\lambda_j}(e^{-i\lambda_j}L)\bar{u}_n\\ &\quad - (2\pi n)^{-1/2} J_n(L)^2 \widetilde{D}_{n\lambda_s}(e^{-i\lambda_j}L;\theta)\bar{u}_n.\end{aligned}$$

By the same argument as the ones used in showing (77) and (82), we obtain

$$E|L^q \widetilde{J}_{n\lambda_j}(e^{-i\lambda_j}L)\bar{u}_n|^2 = O(nj^{-1}), \qquad q = 0, 1, \ldots.$$

In conjunction with $J_n(e^{i\lambda_j}) = O(\log n)$, Minkowski's inequality and (82), it follows that

$$\begin{aligned}&E\sup_\theta |n^{\theta-1/2} j^{1/2-\theta} D_n(e^{i\lambda_j};\theta)[J_n(e^{i\lambda_j}) + J_n(L)]\widetilde{J}_{n\lambda_j}(e^{-i\lambda_j}L)\bar{u}_n|^2\\ &\qquad = O((\log n)^2),\\ &E\sup_\theta |n^{\theta-1/2} j^{1/2-\theta} J_n(L)^2 \widetilde{D}_{n\lambda_j}(e^{-i\lambda_j}L;\theta)\bar{u}_n|^2\\ &\qquad = O((\log n)^6)\end{aligned}$$

for $j = 1, \ldots, m$, giving the stated result. □

LEMMA 5.10. *Let $Q_k$, $k = 0, \ldots, 3$, be any real numbers, $\kappa \in (0, \frac{1}{8})$, and $1/m + m/n \to 0$ as $n \to \infty$. Then there exists $\eta > 0$ not depending on $Q_k$ such that, for sufficiently large $n$:*

(a) $$m^{-1}\sum_{j=[\kappa m]}^{m}|(1-e^{i\lambda_j})^3Q_3+(1-e^{i\lambda_j})^2Q_2+(1-e^{i\lambda_j})Q_1+Q_0|^2 \geq \eta(m^6n^{-6}Q_3^2+m^4n^{-4}Q_2^2+m^2n^{-2}Q_1^2+Q_0^2),$$

(b) $$\begin{aligned} m^{-1}\sum_{j=[\kappa m]}^{m}&|(1-e^{i\lambda_j})^{-1}Q_3+(1-e^{i\lambda_j})^{-2}Q_2 \\ &+(1-e^{i\lambda_j})^{-3}Q_1+(1-e^{i\lambda_j})^{-4}Q_0|^2 \\ &\geq \eta(m^{-2}n^2Q_3^2+m^{-4}n^4Q_2^2+m^{-6}n^6Q_1^2+m^{-8}n^8Q_0^2). \end{aligned}$$

Proof. Define

$$A(\lambda)=(1-e^{i\lambda})^3Q_3+(1-e^{i\lambda})^2Q_2+(1-e^{i\lambda})Q_1+Q_0.$$

Since $1-e^{i\lambda}=-i\lambda+O(\lambda^2)$ as $\lambda\to 0$, we have

$$\begin{aligned} A(\lambda)&=i\lambda^3Q_3-\lambda^2Q_2-i\lambda Q_1+Q_0 \\ &\quad+O(\lambda^4)Q_3+O(\lambda^3)Q_2+O(\lambda^2)Q_1. \end{aligned} \tag{88}$$

Applying $2|a||b|\leq |a|^2+|b|^2$ to the product terms involving the remainder terms, we obtain

$$|A(\lambda)|^2=(\lambda^2Q_2-Q_0)^2+(\lambda^3Q_3-\lambda Q_1)^2+R(\lambda), \tag{89}$$

where $R(\lambda)=O(\lambda^7)Q_3^2+O(\lambda^5)Q_2^2+O(\lambda^3)Q_1^2+O(\lambda)Q_0^2$. First we show that

$$m^{-1}\sum_{j=[\kappa m]}^{m}(\lambda_j^2Q_2-Q_0)^2\geq \eta(m^4n^{-4}Q_2^2+Q_0^2). \tag{90}$$

When $\operatorname{sgn}(Q_2)\neq\operatorname{sgn}(Q_0)$, then (90) holds from Lemma 5.4. When $\operatorname{sgn}(Q_2)=\operatorname{sgn}(Q_0)$, without loss of generality assume $Q_2, Q_0>0$. Note that $\lambda_j^2Q_2$ is an increasing function of $j$. Now suppose $(\lambda_{m/2})^2Q_2-Q_0\geq 0$. Then, since $(\lambda_{3m/4})^2=\frac{9}{4}(\lambda_{m/2})^2$, we have, for $j=3m/4,\dots,m$,

$$\begin{aligned} \lambda_j^2Q_2-Q_0&\geq \tfrac{9}{4}(\lambda_{m/2})^2Q_2-Q_0 \\ &=\tfrac{1}{4}(\lambda_{m/2})^2Q_2+2(\lambda_{m/2})^2Q_2-Q_0 \\ &\geq \tfrac{1}{4}(\lambda_{m/2})^2Q_2+Q_0. \end{aligned}$$

Now suppose $(\lambda_{m/2})^2Q_2-Q_0<0$. Then, since $(\lambda_{m/4})^2=\frac{1}{4}(\lambda_{m/2})^2$, we have, for $j=1,\dots,m/4$,

$$\begin{aligned} \lambda_j^2Q_2-Q_0&\leq \tfrac{1}{4}(\lambda_{m/2})^2Q_2-Q_0 \\ &=-\tfrac{1}{4}(\lambda_{m/2})^2Q_2+[\tfrac{1}{2}(\lambda_{m/2})^2Q_2-Q_0] \\ &\leq -\tfrac{1}{4}(\lambda_{m/2})^2Q_2-\tfrac{1}{2}Q_0. \end{aligned}$$

Therefore, either for $j = 1, \dots, m/4$ or for $j = 3m/4, \dots, m$, we have

$$|\lambda_j^2 Q_2 - Q_0| \geq \tfrac{1}{4}(\lambda_{m/2})^2 Q_2 + \tfrac{1}{2} Q_0 \tag{91}$$

and (90) follows immediately. The same argument gives, if $\operatorname{sgn}(Q_3) = \operatorname{sgn}(Q_1)$,

$$|\lambda_j^3 Q_3 - \lambda_j Q_1| \geq \lambda_j \{\tfrac{1}{4}(\lambda_{m/2})^2 |Q_3| + \tfrac{1}{2}|Q_1|\}, \tag{92}$$

either for $j = 1, \dots, m/4$ or for $j = 3m/4, \dots, m$, and it follows from (91) and (92) that

$$m^{-1} \sum_{j=[\kappa m]}^{m} (\lambda_j^3 Q_3 - \lambda_j Q_1)^2 \geq \eta(m^6 n^{-6} Q_3^2 + m^2 n^{-2} Q_1^2).$$

For $R(\lambda)$ in (89), it follows from Lemma 5.4 that

$$m^{-1} \sum_{j=[\kappa m]}^{m} R(\lambda_j)$$
$$= O(m^7 n^{-7}) Q_3^2 + O(m^5 n^{-5}) Q_2^2 + O(m^3 n^{-3}) Q_1^2 + O(m n^{-1}) Q_0^2,$$

and part (a) follows. For part (b), rewrite the term inside the summation as

$$|(1 - e^{i\lambda_j})^{-4} A(\lambda_j)|^2 = |\lambda_j^{-4}(1 + O(\lambda_j)) A(\lambda_j)|^2.$$

Applying (88) and the following argument with (91) and (92) gives part (b). □

**Acknowledgments.** The authors thank an Associate Editor and three referees for helpful comments and advice that led to a substantial revision of the original version of this paper. K. Shimotsu thanks the Cowles Foundation for hospitality during his stay from January 2002 to August 2003. Simulations were performed in MATLAB.

Department of Economics
Queen's University
Kingston, Ontario
Canada K7L 3N6
e-mail: shimotsu@qed.econ.queensu.ca

Cowles Foundation for
Research in Economics
Yale University
P.O. Box 208281
New Haven, Connecticut 06520-8281
USA
e-mail: peter.phillips@yale.edu